\documentclass[11pt]{amsart}
\usepackage{amscd}
\usepackage{pinlabel}
\usepackage{graphicx}

\newtheorem{thm}{Theorem}[section]
\newtheorem{cnj}[thm]{Conjecture}
\newtheorem{cor}[thm]{Corollary}
\newtheorem{lma}[thm]{Lemma}

\theoremstyle{definition}
\newtheorem{dfn}[thm]{Definition}

\newenvironment{pf}{\begin{proof}}{\end{proof}}

\numberwithin{equation}{section}

\newcommand{\R}{{\mathbb{R}}}
\newcommand{\C}{{\mathbb{C}}}

\newcommand{\Z}{{\mathbb{Z}}}

\newcommand{\V}{{\mathbf{V}}}

\newcommand{\Aa}{{\mathcal{A}}}

\newcommand{\M}{{\mathcal{M}}}
\newcommand{\la}{\langle}
\newcommand{\ra}{\rangle}
\newcommand{\pa}{\partial}
\newcommand{\id}{\operatorname{id}}
\newcommand{\img}{\operatorname{im}}
\newcommand{\krn}{\operatorname{ker}}

\newcommand{\ev}{\operatorname{ev}}

\newcommand{\Hom}{\operatorname{Hom}}
\newcommand{\action}{{\mathfrak a}}

\title[Rational SFT, linearized Legendrian contact homology\dots]
{Rational SFT, linearized Legendrian contact homology, and Lagrangian Floer cohomology}
\author{Tobias Ekholm}
\address{Department of mathematics, Uppsala University, Box 480, 751 06 Uppsala, Sweden}
\email{tobias{\@@}math.uu.se}

\begin{document}
\thanks{The author acknowledges support from the G{\"o}ran Gustafsson Foundation for Research in Natural Sciences and Medicine.}
\maketitle

\begin{center}
{\em This paper is dedicated to Oleg Viro on the occasion of his 60th birthday.}
\end{center}

\begin{abstract}
We relate the version of rational Symplectic Field Theory for exact Lagrangian cobordisms introduced in \cite{E-JEMS} with linearized Legendrian contact homology. More precisely, if $L\subset X$ is an exact Lagrangian submanifold of an exact symplectic manifold with convex end $\Lambda\subset Y$, where $Y$ is a contact manifold and $\Lambda$ is a Legendrian submanifold, and if $L$ has empty concave end, then the linearized Legendrian contact cohomology of $\Lambda$, linearized with respect to the augmentation induced by $L$, equals the rational SFT of $(X,L)$. Following ideas of P.~Seidel \cite{S-private}, this equality in combination with a version of Lagrangian Floer cohomology of $L$ leads us to a conjectural exact sequence which in particular implies that if $X=\C^{n}$ then the linearized Legendrian contact cohomology of $\Lambda\subset S^{2n-1}$ is isomorphic to the singular homology of $L$. We outline a proof of the conjecture and show how to interpret the duality exact sequence for linearized contact homology of \cite{EESa-duality} in terms of the resulting isomorphism.
\end{abstract}

\section{Introduction}\label{S:intr}
Let $Y$ be a contact $(2n-1)$-manifold with contact $1$-form $\lambda$ (i.e., $\lambda\wedge (d\lambda)^{n-1}$ is a volume form on $Y$). The {\em Reeb vector field} $R_\lambda$ of $\lambda$ is the unique vector field which satisfies $\lambda(R_\lambda)=1$ and $d\lambda(R_\lambda,\cdot)=0$. The {\em symplectization} of $Y$ is the symplectic manifold $Y\times\R$ with symplectic form $d(e^t\lambda)$ where $t$ is a coordinate in the $\R$-factor. A {\em symplectic manifold with cylindrical ends} is a symplectic $2n$-manifold $X$ which contains a compact subset $K$ such that $X-K$ is symplectomorphic to a disjoint union of two half-symplectizations $Y^+\times\R_+\cup Y^-\times\R_-$, for some contact $(2n-1)$-manifolds $Y^\pm$, where $\R_+=[0,\infty)$ and $\R_-=(-\infty,0]$. We call $Y^+\times\R_+$  and of $Y^-\times\R_-$ the {\em positive- and negative ends} of $X$, respectively, and $Y^+$ and of $Y^-$, {\em $(+\infty)$- and $(-\infty)$-boundaries} of $X$, respectively.

The relative counterpart of a symplectic manifold with cylindrical ends is a pair $(X,L)$ of a symplectic $2n$-manifold $X$ with a Lagrangian $n$-submanifold $L\subset X$ (i.e., the restriction of the symplectic form in $X$ to any tangent space of $L$ vanishes) such that outside a compact subset, $(X,L)$ is symplectomorphic to the disjoint union of $\bigl(Y^+\times\R_+,\Lambda^+\times\R_+\bigr)$ and $\bigl(Y^-\times\R_-,\Lambda^-\times\R_-\bigr)$, where $\Lambda^\pm\subset Y^{\pm}$ are Legendrian $(n-1)$-submanifolds (i.e., $\Lambda^\pm$ are everywhere tangent to the kernels of the contact forms on $Y^\pm$). If the symplectic manifold $X$ is exact (i.e., if the symplectic form $\omega$ on $X$ satisfies $\omega=d\beta$ for some $1$-form $\beta$) and if the Lagrangian submanifold $L$ is exact as well (i.e., if the restriction $\beta|L$ satisfies $\beta|L=df$ for some function $f$), then we call the pair $(X,L)$ of exact manifolds an {\em exact cobordism}. We assume throughout the paper that $X$ is simply connected, that the first Chern class of $TX$, viewed as a complex bundle using any almost complex structure compatible with the symplectic form on $X$, is trivial and that the Maslov class of $L$ is trivial as well. (These assumptions are made in order to have well defined gradings in contact homology algebras over $\Z_2$. In more general cases one would work with contact homology algebras with suitable Novikov coefficients in order to have appropriate gradings.)

In \cite{E-JEMS} a version of rational Symplectic Field Theory (SFT) (see \cite{EGH} for a general description of SFT) for exact cobordisms with good ends, see Section \ref{S:relSFT}, was developed. (The additional condition that ends be good allows us to disregard Reeb orbits in the ends when setting up the theory. Standard contact spheres as well as $1$-jet spaces with their standard contact structures are good.) It associates to  an exact cobordism $(X,L)$, where $L$ has $k$ components, a $\Z$-graded filtered $\Z_2$-vector space $\V(X,L)$, with $k$ filtration levels, and with a filtration preserving differential $d^{f}\colon \V(X,L)\to\V(X,L)$. Elements in $\V(X,L)$ are formal sums of admissible formal disks in which the number of summands with $(+)$-action below any given number is finite, see Section \ref{S:relSFT} for definitions of these notions. The differential increases $(+)$-action and hence if $\V_{[\alpha]}(X,L)$ denotes $\V(X,L)$ divided out by the subcomplex of all formal sums in which all disks have $(+)$-action larger than $\alpha$, then the differential induces filtration preserving differentials $d^{f}_{\alpha}\colon \V_{[\alpha]}(X,L)\to\V_{[\alpha]}(X,L)$ with associated spectral sequences $\bigl\{E^{p,q}_{r;[\alpha]}(X,L)\bigr\}_{r=1}^{k}$. The projection maps $\pi^{\alpha}_{\beta}\colon \V_{[\alpha]}(X,L)\to\V_{[\beta]}(X,L)$, $\alpha>\beta$ give an inverse system of chain maps. The limit $E_{r}^{\ast}(X,L)=\underleftarrow{\lim}_{\alpha} E_{r;[\alpha]}^{\ast}(X,L)$ is invariant under deformations of $(X,L)$ through exact cobordisms with good ends and in particular under deformations  of $L$ through exact Lagrangian submanifolds with cylindrical ends, see Theorem \ref{t:rSFTinv}.

In this paper we will use only the simplest version of the theory just described which is as follows. Let $(X,L)$ be an exact cobordism such that $L$ is connected and without negative end, i.e., $\Lambda^{-}=\emptyset$. In this case, admissible formal disks have only one positive puncture and we identify (a quotient of) $\V(X,L)$ with the $\Z_2$-vector space of formal sums of Reeb chords of $\Lambda=\Lambda^{+}$. Furthermore our assumptions on $\pi_1(X)$ and vanishing of $c_1(TX)$ and of the Maslov class of $L$ imply that the grading of a formal disks depends only on the Reeb chord at its positive puncture. We let $|c|$ denote the grading of a chord $c\in\V(X,L)$. Rational SFT then provides a differential
\[
d^{f}\colon \V(X,L)\to\V(X,L)
\]
with $|d^{f}(c)|=|c|+1$, and
which increases action in the sense that if $\action(c)$ denotes the action of the Reeb chord $c$ and if the Reeb chord $b$ appears with non-zero coefficient in $\pa c$ then $\action(b)>\action(c)$. Furthermore, since $L$ is connected the spectral sequences have only one level and
\[
E_{1}^{\ast}(X,L)=\underleftarrow\lim_{\alpha}\left(\krn d^{f}_{\alpha}/\img d^{f}_\alpha\right).
\]

Our first result relates $E_{1}^{\ast}(X,L)$ to linearized Legendrian contact cohomology, see Section \ref{S:LCH}. Legendrian contact homology was introduced in \cite{EGH, Ch}. It was worked out in detail in special cases including $1$-jet spaces in \cite{EES-JDGanalysis,EES-orientations,EES-PxR}. From the point of view of Legendrian contact homology an exact cobordism $(X,L)$ with good ends induces a chain map from the contact homology algebra of  $(Y^{+},\Lambda^{+})$ to that of $(Y^{-},\Lambda^{-})$, in particular if $\Lambda^{-}=\emptyset$ then the latter equals the ground field $\Z_{2}$ with the trivial differential. Such a chain map $\epsilon$ is called an {\em augmentation} and gives rise to a linearization of the contact homology algebra of $\Lambda^{+}$. That is, it endows the chain complex $Q(\Lambda)$ generated by Reeb chords of $\Lambda$ with a differential $\pa^{\epsilon}$. The resulting homology is called {\em $\epsilon$-linearized contact homology} and denoted $LCH_{\ast}(Y,\Lambda;\epsilon)$. We let $LCH^{\ast}(Y;\Lambda;\epsilon)$ be the homology of the dual complex $Q'(\Lambda)=\Hom(Q(\Lambda);\Z_2)$ and call it the {\em $\epsilon$-linearized contact cohomology} of $\Lambda$.

We say that $(X,L)$ satisfies a {\em monotonicity condition} if there are constants $C_0$ and $C_1>0$ such that for any Reeb chord $c$ of $\Lambda\subset Y$, $|c|> C_1\action(c) +C_0$. Note that if $Y$ is a $1$-jet space or the sphere endowed with a generic small perturbation of the standard contact form and if $\Lambda$ is in general position with respect to the Reeb flow then $(X,L)$ satisfies a monotonicity condition.

\begin{thm}\label{t:iso}
Let $(X,L)$ be an exact cobordism with good ends. Let $(Y,\Lambda)$ denote the positive end of $(X,L)$ and assume that the $(-\infty)$-boundary of $L$ is empty. Let $\epsilon$ denote the augmentation on the contact homology algebra of $\Lambda$ induced by $L$. Then the natural map $Q'(\Lambda)\to \V(X,L)$, which takes an element in $Q'(\Lambda)$ thought of as a formal sum of co-vectors dual to Reeb chords in $Q(\Lambda)$ to the corresponding formal sum of Reeb chords in $\V(X,L)$, is a chain map. Furthermore if $(X,L)$ satisfies a monotonicity condition then the corresponding map on homology
\[
LCH^{\ast}(Y,\Lambda;\epsilon)\to E^{\ast}_{1}(X,L),
\]
is an isomorphism.
\end{thm}

Theorem \ref{t:iso} is proved in Section \ref{s:pft:iso}. We point out that when $(X,L)$ satisfies a monotonicity condition it follows from this result that $LCH_{\ast}(Y,\Lambda;\epsilon)$ depends only on the symplectic topology of $(X,L)$.

We next consider two exact cobordisms $(X,L_0)$ and $(X,L_1)$ with good ends and with the following properties: both $L_0$ and $L_1$ have empty $(-\infty)$-boundaries, if $\Lambda_j$ denotes the $(+\infty)$-boundary of $L_j$ then $\Lambda_0\cap \Lambda_1=\emptyset$, $L_0$ and $L_1$ intersect transversely, and the Reeb flow of $\Lambda_0$ along a Reeb chord connecting $\Lambda_0$ to $\Lambda_1$ is transverse to $\Lambda_1$ at its endpoint. For such pairs of exact cobordisms we  define Lagrangian Floer cohomology $HF^{\ast}(X;L_0,L_1)$ as an inverse limit of the homologies $HF^{\ast}_{[\alpha]}(X;L_0,L_1)$ of chain complexes $C_{[\alpha]}(X;L_0,L_1)$ generated by Reeb chords between $\Lambda_0$ and $\Lambda_1$ of action at most $\alpha$ and by points in $L_0\cap L_1$. This Floer cohomology has a relative $\Z$-grading and is invariant under exact deformations of $L_1$.

Consider an exact cobordism $(X,L)$ where $L$ has empty $(-\infty)$-boundary and $(+\infty)$-boundary $\Lambda$. Let $L'$ be a slight push off of $L$ which is an extension of a small push off $\Lambda'$ of $\Lambda$ along the Reeb vector field.
\begin{cnj}\label{cnj:Seidel}
For any $\alpha>0$, there is a long exact sequence
\begin{equation}\label{e:longexact}
\begin{CD}
\cdots @>{\delta_{\alpha;L,L'}}>> E_{1;[\alpha]}^{\ast}(X,L) @>>> HF_{[\alpha]}^{\ast}(X;L,L')
@>>> H_{n-\ast}(L)\\
@>{\delta_{\alpha;L,L'}}>> E_{1;[\alpha]}^{\ast+1}(X,L) @>>> HF_{[\alpha]}^{\ast+1}(X;L,L') @>>> H_{n-\ast-1}(L)\\\
@>{\delta_{\alpha;L,L'}}>> \cdots,
\end{CD}
\end{equation}
where $H_\ast(L)$ is the ordinary homology of $L$ with $\Z_2$-coefficients. It follows in particular, that if $X=\C^{n}$ or $X=J^{1}(\R^{n-1})\times\R$ then $HF^{\ast}(X;L,L')=0$ and the map $\delta_{L,L'}\colon H_{n-\ast+1}(L)\to E_1^{\ast}(X,L)\approx LCH^{\ast}(Y,\Lambda;\epsilon)$, induced by the maps $\delta_{\alpha;L,L'}$, is an isomorphism.
\end{cnj}

The author learned about the isomorphism above, between linearized contact homology of a Legendrian submanifold with a Lagrangian filling and the ordinary homology of the filling, from P.~Seidel \cite{S-private} who explained it using an exact sequence in wrapped Floer homology \cite{AS, FSS} similar to \eqref{e:longexact}. Borrowing Seidel's argument, we outline a proof of Conjecture \ref{cnj:Seidel} in Section \ref{s:pfcnj:Seidel} in which the Lagrangian Floer cohomology $HF(X;L,L')$ plays the role of wrapped Floer homology.

In \cite{EESa-duality} a duality exact sequence for linearized contact homology of a Legendrian submanifold $\Lambda\subset Y$, where $Y=P\times\R$ for some exact symplectic manifold $P$ and where the projection of $\Lambda$ into $P$ is displaceable, was found. In what follows we restrict attention to the case $Y=J^{1}(\R^{n-1})$. Then every compact Legendrian submanifold has displaceable projection and the duality exact sequence is the following, where $\epsilon$ denotes any augmentation and where we suppress the ambient manifold $Y=J^{1}(\R^{n-1})$ from the notation,
\begin{equation}\label{e:seq}
\begin{CD}
\cdots @>{\rho}>> H_{k+1}(\Lambda) @>{\sigma}>>
    LCH^{(n-1)-k-1}(\Lambda;\epsilon) @>{\theta}>> LCH_k(\Lambda;\epsilon)\\
    @>{\rho}>> H_k(\Lambda) @>{\sigma}>>
    LCH^{(n-1)-k}(\Lambda;\epsilon) @>{\theta}>> LCH_{k-1}(\Lambda;\epsilon)
\cdots.
\end{CD}
\end{equation}
Here, if $\beta=\rho(\alpha)\in H_k(\Lambda)$, then the Poincar\'e dual
$\gamma\in H_{n-k}(\Lambda)$ of $\beta$ satisfies $\la\sigma(\gamma),\alpha \ra=1$, where $\la\,,\ra$ is the pairing between the homology and cohomology of $Q(\Lambda)$. Furthermore, the maps $\rho$ and $\sigma$ are defined through a count of rigid configurations of holomorphic disks with boundary on $\Lambda$ with a flow line emanating from its boundary, and the map $\theta$ is defined through a count of rigid holomorphic disks with boundary on $\Lambda$ with two positive punctures.

In $J^{1}(\R^{n-1})$, a generic Legendrian submanifold has finitely many Reeb chords. Furthermore, if $L$ is an exact Lagrangian cobordism in the symplectization $J^{1}(\R^{n-1})\times\R$ with empty $(-\infty)$-boundary then $L$ is displaceable. Hence, both Theorem \ref{t:iso} and Conjecture \ref{cnj:Seidel} give isomorphisms. Combining \eqref{e:longexact} and \eqref{e:seq} leads to the following.
\begin{cor}\label{c:bigdiagram}
Let $L$ be an exact Lagrangian cobordism in $J^{1}(\R^{n-1})\times\R$ with empty $(-\infty)$-boundary and with $(+\infty)$-boundary $\Lambda$ and let $\epsilon$ denote the augmentation on the contact homology algebra of $\Lambda$ induced by $L$. Then the following diagram with exact rows commutes and all vertical maps are isomorphisms
\[
\cdots\,
\begin{CD}
H_{k+1}(\Lambda) @>>>
H_{k+1}(L) @>>>
H_{k+1}(L,\Lambda) @>>> H_{k}(\Lambda)\\
@V{\id}VV @V{\delta_{L,L'}}VV @VV{H^{-1}\circ \delta_{L,L'}'}V @VV{\id}V\\
H_{k+1}(\Lambda) @>{\sigma}>>
LCH^{n-k-2}(\Lambda;\epsilon) @>>>
LCH_k(\Lambda;\epsilon) @>{\rho}>> H_k(\Lambda)\\
\end{CD}
\cdots
\]
Here the top row is the long exact homology sequence of $(L,\Lambda)$, the bottom row is the duality exact sequence, the map $\delta_{L,L'}$ is the map in Conjecture \ref{cnj:Seidel}, the map $\delta'_{L,L'}$ is analogous to $\delta_{L,L'}$, and the map $H$ counts disks in the symplectization with boundary on $\Lambda$ and with two positive punctures, see Section \ref{s:pfc:bigdiagram} for details.
\end{cor}
The proof of Corollary \ref{c:bigdiagram} is discussed in Section \ref{s:pfc:bigdiagram}.

\section{A brief sketch of relative SFT of Lagrangian cobordisms}\label{S:relSFT}
Although we will only use the simplest version of relative SFT introduced in \cite{E-JEMS} in this paper, we give a brief introduction to the full theory for two reasons. First, it is reasonable to expect that this theory is related to product structures on linearized contact homology, see \cite{CEKSW}, in much the same way as the simplest version of the theory appears in Conjecture \ref{cnj:Seidel}. Second, some of the moduli spaces of holomorphic disks that we will make use of are analogous to those needed for more involved versions of the theory.

\subsection{Formal- and admissible disks}\label{s:formdi}
In order to describe relative rational SFT we introduce the following notation. Let $(X,L)$ be an exact cobordism with ends $(Y^{\pm}\times\R_{\pm},\Lambda^\pm\times\R_\pm)$. Write $(\bar X,\bar L)$ for a compact part of $(X,L)$ obtained by cutting the infinite parts of the cylindrical ends off at some $|t|=T>0$. We will sometimes think of Reeb chords of $\Lambda^\pm$ in the $(\pm\infty)$-boundary as lying in $\pa \bar X$ with endpoints on $\pa\bar L$. A {\em formal disk} of $(X,L)$ is a homotopy class of maps of the $2$-disk $D$, with $m$ marked disjoint closed subintervals in $\pa D$, into $\bar X$, where the $m$ marked intervals are required to map in an orientation preserving (reversing) manner to Reeb chords of $\pa \bar L$ in the $(+\infty)$-boundary (in the $(-\infty)$-boundary) and where remaining parts of the boundary $\pa D$ map to $\bar L$.

If $L^{b}$ and $L^{a}$ are exact Lagrangian cobordisms in $X^{b}$ and $X^{a}$, respectively such that a component $(Y,\Lambda)$ of the $(-\infty)$-boundary of $(X^{a},L^{a})$ agrees with a component of the  $(+\infty)$-boundary of $(X^{b},L^{b})$ then these cobordisms can be joined to an exact cobordism $L^{ba}$ in $X^{ba}$, where $(X^{ba},L^{ba})$ is obtained by gluing the positive end $(Y,\Lambda)$ of $(X^{b},L^{b})$ to the corresponding negative end of $(X^{a},L^{a})$. Furthermore if $v^{b}$ and $v^{a}$ are collections of formal disks of $(X^{b},L^{b})$ and $(X^{a},L^{a})$, respectively, then we can construct formal disks in $L^{ba}$ in the following way: start with a disk $v^{a}_1$ from $v^{a}$, let $c_1,\dots,c_{r_1}$ denote the Reeb chords at its negative punctures. Attach positive punctures of disks $v^{b}_{1;1},\dots, v^{b}_{1;r_1}$ in $v^{b}$ mapping to the Reeb chords $c_1,\dots,c_r$ to the corresponding negative punctures of the disk $v^{a}_1$. This gives a disk $v^{ba}_{1}$ with some positive punctures mapping to chords $c_1,\dots,c_{r_2}$ of $\Lambda$. Attach negative punctures of disk $v^{a}_{2;1},\dots,v^{a}_{2;r_2}$ to $v^{ba}_1$ at $c_1,\dots,c_{r_2}$. This gives a disk $v^{ba}_2$ with some negative punctures mapping to Reeb chords in $\Lambda$. Continue this process until there are no punctures mapping to $\Lambda$. We call the resulting disk a formal disk in $L^{ba}$ {\em with factors from $v^{a}$ and $v^{b}$}.

Assume that the set of connected components of $L$ has been subdivided into subsets $L_j$ so that $L$ is a disjoint union $L=L_1\cup\dots\cup L_k$ where each $L_j$ is a collection of connected components of $L$. We call $L_1,\dots,L_k$ the {\em pieces} of $L$. With respect to such a subdivision, Reeb chords fall into two classes: {\em pure} with both endpoints on the same piece and {\em mixed} with endpoints on distinct pieces.

A formal disk represented by a map $u\colon D\to \bar X$ is {\em admissible} if for any arc $\alpha$ in $D$ which connects two unmarked segments in $\pa D$ that are mapped to the same piece by $u$, all marked segments on the boundary of one of the components of $D-\alpha$ map to pure Reeb chords in the $(-\infty)$-boundary.

\subsection{Holomorphic disks}\label{s:hodi}
Let $(X,L)$ be an exact cobordism. Fix an almost complex structure $J$ on $X$ which is adjusted to its symplectic form. Let $S$ be a punctured Riemann surface with complex structure $j$ and with boundary $\pa S$. A $J$-holomorphic curve with boundary on $L$ is a map $u\colon S\to X$ such that
\[
du+J\circ du\circ j=0,
\]
and such that $u(\pa S)\subset L$. For details on holomorphic curves in this setting we refer to \cite[Appendix B]{E-JEMS} and references therein. Here we summarize the main properties we will use.
By definition, an adjusted almost complex structure $J$ is invariant under $\R$-translations in the ends of $X$ and pairs the Reeb vector field in the $(\pm\infty)$-boundary with the symplectization direction. Consequently, strips which are cylinders over Reeb chords as well as cylinders over Reeb orbits are $J$-holomorphic. Furthermore, any $J$-holomorphic disk of finite energy is asymptotic to such Reeb chord strips at its boundary punctures and to Reeb orbit cylinders at interior punctures, see \cite[Section B.1]{E-JEMS}. We say that a puncture of a $J$-holomorphic disk is positive (negative) if the disk is asymptotic to a Reeb chord strip (Reeb orbit cylinder) in the positive (negative) end of $(X,L)$. Note that exactness of $(X,L)$ and the fact that the symplectic form is positive on $J$-complex tangent lines imply that any $J$-holomorphic curve has at least one positive puncture.

These results on asymptotics imply that any $J$-holomorphic disk in $X$ with boundary on $L$ determines a formal disk. Let $\M(v)$ denote the moduli space of $J$-holomorphic disks with associated formal disk equal to $v$. The formal dimension of $\M(v)$ is determined by the Fredholm index of the linearized $\bar\pa_{J}$-operator along a representative of $v$, see \cite[Section 3.1]{E-JEMS}.

A sequence of $J$-holomorphic disks with boundary on $L$ may converge to a broken disk of two components which intersects at a boundary point. We will refer to this phenomenon as boundary bubbling. However, if all elements in the sequence have only one positive puncture then boundary bubbling is impossible by exactness: each component in the limit curve must have at least one positive puncture. The reason for using admissible disks to set up relative SFT is the following:
in a sequence of holomorphic disks with corresponding formal disks admissible, boundary bubbling is impossible for topological reasons. As a consequence, if $v$ is a formal disk then the boundary of $\M(v)$ consists of several level $J$-holomorphic disks and spheres joined at Reeb chords or at Reeb orbits, see \cite{BEHWZ}.

Recall from Section \ref{S:intr} that we require the ends of our exact cobordisms $(X,L)$ to be good. The precise formulation of this condition is as follows. If $\gamma^{+}$ ($\gamma^{-}$) is a Reeb orbit in the $(+\infty)$-boundary $Y^{+}$ (in the $(-\infty)$-boundary $Y^{-}$) of $X$ then the formal dimension of any moduli space of holomorphic spheres in $X$ (in $Y^{-}\times\R$) with positive puncture at $\gamma^{+}$ (at $\gamma^{-}$) is $\ge 2$. Together with transversality arguments these conditions guarantee that broken curves in the boundary of $\M(v)$, where $v$ is an admissible formal disk cannot contain any spheres, if $\dim(\M(v))\le 1$ ($\dim(\M(v))=2$ if $(X,L)$ is a trivial cobordism), see \cite[Lemma B.6]{E-JEMS}. In particular, in the boundary of $\M(v)$, where $\dim(\M(v))$ satisfies these dimensional constraints and where $v$ is admissible, there can be only two level curves, all pieces of which are admissible disks, see \cite[Lemma 2.5]{E-JEMS}.

Under our additional assumptions ($\pi_1(X)$ trivial, first Chern class of $X$ and Maslov class of $L$ vanish) the grading of a formal disk depends only on the Reeb chords at its punctures. For later reference, we describe this more precisely in the case when $L$ is connected and when its $(-\infty)$-boundary is empty. Let $(Y,\Lambda)$ denote the $(+\infty)$-boundary of $(X,L)$. If $c$ is a Reeb chord of $\Lambda\subset Y$ then let $\gamma$ be any path in $L$ joining its endpoints. Since $X$ is simply connected $\gamma\cup c$ bounds a disk $\Gamma\colon D\to X$. Fix a trivialization of $TX$ along $\Gamma$ such that the linearized Reeb flow along $c$ is represented by the identity transformation with respect to this trivialization. Then the tangent space $T_{s}(\Lambda\times\R)$ at the initial point $s$ of $c$ is transported to a subspace $V_s\times\R$ in the tangent space $T_eX$ at the finial point $e$ of $c$ where $V_s$ is transverse to $T_e\Lambda$ in the contact hyperplane $\xi_e$ at $e$. Let $R$ denote a negative rotation along the complex angle taking $V_{s}$ to $T_{e}\Lambda$ in $\xi_e$, see \cite[Section 3.1]{E-JEMS}, then the Lagrangian tangent planes of $L$ along $\gamma$ capped off with $R$ form a loop $\Delta_\Gamma$ of Lagrangian subspaces in $\C^{n}$ with respect to the trivialization and if $\M(c)$ denotes the moduli space of holomorphic disks in $X$ with one positive boundary puncture at which they are asymptotic to the Reeb chord strip of the Reeb chord $c$ then
\[
\dim(\M(c))=n-3+\mu(\Delta_\Gamma)+1,
\]
see \cite[p.~655]{E-JEMS}. To see that this is independent of $\Gamma$, note that the difference of two trivializations along the disks is measured by $c_1(TX)$. To see that it is independent of the path $\gamma$, note that the difference in the dimension formula corresponding to two different paths $\gamma$ and $\gamma'$ is measured by the Maslov class of $L$ evaluated on the loop $\gamma\cup -\gamma'$. Define
\begin{equation}\label{e:gradingSFT}
|c|=\dim(\M(c)).
\end{equation}

If $a$ and $b_1,\dots,b_m$ are Reeb chords of $\Lambda$ and if $\M(a;b_1,\dots,b_m)$ denotes the moduli space of holomorphic disks in $Y\times\R$ with boundary on $\Lambda\times\R$ then additivity of the index gives
\begin{equation}\label{e:gradingCH}
\dim(\M(a;b_1,\dots,b_k))=|a|-\sum_j|b_j|.
\end{equation}

\subsection{Hamiltonian- and potential vectors and differentials}
Let $(X,L)$ be an exact cobordism and let $v$ be a formal disk of $(X,L)$. Define the $(+)$-action of $v$ as the sum of the actions
\[
\action(c)=\int_c\lambda^{+},
\]
over the Reeb chords $c$ at its positive punctures. Here $\lambda^{+}$ is the contact form in the $(+\infty)$-boundary $Y^{+}$ of $X$. Note that for generic Legendrian $(+\infty)$-boundary, $\Lambda^{+}\subset Y^{+}$, the set of actions of Reeb chords is a discrete subset of $\R$. Let $\V(X,L)$ denote the $\Z$-graded vector space over $\Z_2$ with elements which are formal sums of admissible formal disks which contain only a finite number of summands below any given $(+)$-action. The grading on $\V(X,L)$ is the following: the degree of a formal disk $v$ is the formal dimension of the moduli space $\M(v)$ of $J$-holomorphic disks homotopic to the formal disk. We use the natural filtration
\[
0\subset F^k\V(X,L)\subset\dots\subset F^2\V(X,L)\subset F^1\V(X,L)=\V(X,L)
\]
of $\V(X,L)$, where $k$ is the number of pieces of $L$ and where the filtration level is determined by the number of positive punctures. (It is straightforward to check that an admissible formal disk has at most $k$ Reeb chords at the positive end).

We will define a differential $d^{f}\colon \V(X,L)\to\V(X,L)$ which respects this filtration using $1$-dimensional moduli spaces of holomorphic disks. To this end, fix an almost complex structure $J$ on $X$ which is compatible with the symplectic form and adjusted to $d(e^{t}\lambda^{\pm})$ in the ends, where $\lambda^{\pm}$ is the contact form in $(\pm\infty)$-boundary. Assume that $J$ is generic with respect to $0$- and $1$-dimensional moduli spaces of holomorphic disks, see \cite[Lemma B.8]{E-JEMS}. Since $J$ is invariant under translations in the ends, $\R$ acts on moduli spaces $\M(u)$, where $u$ is a formal disk of $(Y^{\pm}\times\R,\Lambda^{\pm}\times\R)$. In this case we define the reduced moduli spaces as $\widehat{\M}(u)=\M(u)/\R$. Let $h^{\pm}\in\V(Y^{\pm}\times\R,\Lambda^{\pm}\times\R)$ denote the vector of admissible formal disks in $Y^{\pm}\times\R$ with boundary on $\Lambda^{\pm}\times\R$ represented by $J$-holomorphic disks:
\begin{equation}\label{e:hamvect}
h^{\pm}=\sum_{\dim(\widehat{\M}(v))=0} |\widehat{\M}(v)|\,v\in\V(\Lambda^{\pm}\times\R),
\end{equation}
where the sum ranges over all formal disks of $\Lambda^{\pm}\times\R$ and where $|\widehat{\M}|$ denotes the $(\text{mod }{2})$-number of points in the compact $0$-manifold $\widehat{\M}$. We call $h^{+}$ and $h^{-}$ the {\em Hamiltonian vectors} of the positive and negative ends, respectively. Similarly, let $f$ denote the generating function of rigid disks in the cobordism:
\begin{equation}\label{e:potvect}
f=\sum_{\dim(\M(v))=0} |\M(v)|\, v\in\V(X,L),
\end{equation}
where the sum ranges over all formal disks of $(X,L)$. We call $f$ the {\em potential vector of $(X,L)$}.

We view elements $w$ in $\V(Y^{\pm}\times\R,\Lambda^{\pm}\times\R)$ and $\V(X,L)$ as sets of admissible formal disks, where the set consists of those formal disks which appear with non-zero coefficient in $w$. Define the differential $d^{f}\colon \V(X,L)\to\V(X,L)$ as the linear map such that if $v$ is an admissible formal disk (a generator of $\V(X,L)$) then $d^{f}(v)$ is the sum of all admissible formal disks obtained in the following way.
\begin{itemize}
\item[$({\rm i})$] Attach a positive puncture of $v$ to a negative puncture of a $h^{+}$-disk, or
\item[$({\rm i}')$] attach a negative puncture of $v$ to a positive puncture of a $h^{-}$-disk.
\item[$({\rm ii})$] If the first step was $({\rm i})$ then attach $f$-disks at remaining negative punctures of the $h^{+}$-disk, or
\item[$({\rm ii}')$] if the first step was $({\rm i}')$ then attach $f$-disks at remaining positive punctures of the $h^{-}$-disk.
\end{itemize}

The fact that this is a differential is a consequence of the product structure of the boundary of the moduli space mentioned above in the case of $1$-dimensional moduli spaces, see \cite[Lemma 3.7]{E-JEMS}. Furthermore, the differential increases grading by $1$ and respects the filtration since any disk in $h^{\pm}$ or in $f$ has at least one positive puncture.

\subsection{The rational admissible SFT spectral sequence}
Fix $\alpha>0$. If $\V_{[\alpha+]}(X,L)\subset \V(X,L)$ denotes the subspace of formal sums of formal disks with $(+)$-action at least $\alpha$ then since holomorphic disks have positive symplectic area it follows that $d^{f}\bigl(\V_{[\alpha+]}(X,L)\bigr)\subset\V_{[\alpha+]}(X,L)$. If $\V_{[\alpha]}(X,L)=\V(X,L)/\V_{[\alpha+]}(X,L)$ then $\V_{[\alpha]}(X,L)$ is isomorphic to the vector space generated by formal disks of $(+)$-action $<\alpha$ and there is a short exact sequence of chain complexes
\[
\begin{CD}
0  @>>> \V_{[\alpha+]}(X,L) @>>> \V(X,L) @>>> \V_{[\alpha]}(X,L) @>>> 0.
\end{CD}
\]
The quotients $\V_{[\alpha]}(X,L)$ form an inverse system
\[
\pi^{\alpha}_{\beta}\colon \V_{[\alpha]}(X,L)\to\V_{[\beta]}(X,L),\quad \alpha>\beta,
\]
of graded chain complexes, where $\pi^{\beta}_{\alpha}$ are the natural projections.
Consequently, the $k$-level spectral sequences corresponding to the filtrations
\[
0\subset F^k\V_{[\alpha]}(X,L)\subset\dots\subset F^2\V_{[\alpha]}(X,L)\subset F^1\V_{[\alpha]}(X,L)=\V_{[\alpha]}(X,L)
\]
which we denote
\[
\bigl\{E^{p,q}_{r;[\alpha]}(X,L)\bigr\}_{r=1}^k
\]
form an inverse system as well and we define
the rational admissible SFT invariant as
\[
\bigl\{E^{p,q}_{r}(X,L)\bigr\}_{r=1}^k=
\underleftarrow{\lim}_{\alpha}\bigl\{E^{p,q}_{r;[\alpha]}(X,L)\bigr\}_{r=1}^k.
\]
This is in general not a spectral sequence but it is under some finiteness conditions. The following result is a consequence of \cite[Theorems 1.1 and 1.2]{E-JEMS}.
\begin{thm}\label{t:rSFTinv}
Let $(X,L)$ be an exact cobordism with a subdivision $L=L_1\cup\dots\cup L_k$ into pieces. Then $\bigl\{E^{p,q}_{r}(X,L)\bigr\}$ does not depend on the choice of adjusted almost complex structure $J$, and is invariant under deformations of $(X,L)$ through exact cobordisms with good ends.
\end{thm}

\begin{pf}
Any such deformation can be subdivided into a compactly supported deformation and a Legendrian isotopy at infinity. The former type of deformations are shown to induce isomorphisms of $\bigl\{E^{p,q}_{r}(X,L)\bigr\}$ in \cite[Theorem 1.1.]{E-JEMS}.

To show that the later type of deformation induces an isomorphism, we note that it gives rise to an invertible exact cobordism, see \cite[Appendix A]{E-JEMS}, and use the same argument as in the proof of \cite[Theorem 1.2.]{E-JEMS} as follows. Let $C_{01}$ be the exact cobordism of the Legendrian isotopy at infinity and let $C_{10}$ be its inverse cobordism. We use the symbol $A\# B$ to denote the result of joining two cobordisms along a common end. Consider first the cobordism
\[
L\# C_{01}\# C_{10}.
\]
Since this cobordism can be deformed by a compact deformation to $L$ we find that the composition of the maps $\Phi\colon \V(X,L\# C_{01})\to\V(X,L\# C_{01}\# C_{10})$ and
$\Psi\colon \V(X,L)\to\V(X,L\# C_{01})$ is chain homotopic to identity. Hence $\Psi$ is injective on homology. Consider second the cobordism
\[
L\# C_{01}\# C_{10}\# C_{01}.
\]
Since this cobordism can be deformed to $L\# C_{01}$ we find similarly that there is a map $\Theta$ such that $\Psi\circ\Theta$ is chain homotopic to the identity on $\V(X,L\# C_{01})$ hence $\Psi$ is surjective on homology as well.
\end{pf}

\subsection{A simple version of rational admissible SFT}\label{s:simpleSFT}
As mentioned in Section \ref{S:intr}, in the present paper, we will use the rational admissible spectral sequence in the simplest case: for $(X,L)$ where $L$ has only one component. Since there is only one piece, the spectral sequence has only one level and
\[
E^{1+q}_{1}(X,L)=\underleftarrow{\lim}_{\alpha}\,E_{1;[\alpha]}^{1,q}(X,L)=\underleftarrow{\lim}_{\alpha}\,\krn(d^{f}_{\alpha})/\img(d^{f}_{\alpha})
\]
is the invariant we will compute. To simplify things further we will work not with the chain complex $\V(X,L)$ as described above but with the quotient of it obtained by forgetting the homotopy classes of formal disks. We view this quotient, using our assumption $\pi_1(X)$ trivial, as the space of formal sums of Reeb chords of the $(+\infty)$-boundary $\Lambda$ of $L$. Further, our assumptions $c_1(TX)=0$ and vanishing Maslov class of $L$ implies that the grading descends to the quotient, see \eqref{e:gradingSFT}. For simplicity we keep the notation $\V(X,L)$ and $\V_{[\alpha]}(X,L)$ for the corresponding quotients.

\section{Legendrian contact homology, augmentation, and linearization}\label{S:LCH}
In this section we will define Legendrian contact homology and its linearization. We work in the following setting: $(X,L)$ is an exact cobordism with good ends, the $(-\infty)$-boundary of $L$ is empty and the $(+\infty)$-boundary of $(X,L)$ will be denoted $(Y,\Lambda)$.

Recall that the assumption on good ends allows us to disregard Reeb orbits. Furthermore, our additional assumptions on $(X,L)$, $\pi_1(X)$ trivial and first Chern class and Maslov class trivial, allows us to work with coefficients in $\Z_2$ and still retain grading.

\subsection{Legendrian contact homology}
Assume that $\Lambda\subset Y$ is generic with respect to the Reeb flow on $Y$. If $c$ is a Reeb chord of $\Lambda$, let $|c|\in\Z$ be as in \eqref{e:gradingSFT}.

\begin{dfn}
The DGA of $(Y,\Lambda)$ is the unital algebra $\Aa(Y,\Lambda)$ over $\Z_2$ generated by the Reeb chords of $\Lambda$. The grading of a Reeb chord $c$ is $|c|$.
\end{dfn}

\begin{dfn}\label{d:CHdiff}
The contact homology differential is the map $\pa\colon \Aa(Y,\Lambda)\to\Aa(Y,\Lambda)$ which is linear over $\Z_2$, which satisfies Leibniz rule, and which is defined as follows on generators:
\[
\pa c = \sum_{\dim(\M(c;\overline{b}))=1}
|\widehat{\M}(c;\overline{b})|\,\overline{b},
\]
where $c$ is a Reeb chord and $\overline{b}=b_1\dots b_k$ is a word of Reeb chords. (For notation, see \eqref{e:gradingCH}.)
\end{dfn}

We give a brief explanation of why $\pa$ in Definition \ref{d:CHdiff} is a differential, i.e., why $\pa^{2}=0$. Consider the boundary of the $2$-dimensional moduli space $\M$ (which become $1$-dimensional after the $\R$-action has been divided out) of holomorphic disks with one positive puncture at $a$. As explained in Section \ref{s:hodi}, the boundary of such a moduli space consists of two level curves with all components except two Reeb chord strips. Since these configurations are exactly what is counted by $\pa^{2}c$ and since they correspond to the boundary points of the compact $1$-manifold $\M/\R$, we conclude that $\pa^{2}c=0$.

\begin{dfn}
An {\em augmentation} of $\Aa(Y,\Lambda)$ is a chain map $\epsilon\colon\Aa(Y,\Lambda)\to\Z_{2}$, where $\Z_{2}$ is equipped with the trivial differential.
\end{dfn}

Given an augmentation $\epsilon$, define the algebra isomorphism $E_\epsilon\colon\Aa(Y,\Lambda)\to\Aa(Y,\Lambda)$ by letting
\[
E_\epsilon(c)=c+\epsilon(c),
\]
for each generator $c$. Consider the word length filtration of $\Aa(Y,\Lambda)$,
\[
\Aa(Y,\Lambda)=\Aa_{0}(Y,\Lambda)\supset\Aa_1(Y,\Lambda)\supset\Aa_{2}(Y,\Lambda)\supset\dots.
\]
The differential $\pa^{\epsilon}=E_\epsilon\circ\pa\circ E_{\epsilon}^{-1}\colon\Aa(Y,\Lambda)\to\Aa(Y,\Lambda)$ respects this filtration: $\pa^{\epsilon}(\Aa_j(Y,\Lambda))\subset\Aa_j(Y,\Lambda)$. In particular, we obtain the {\em $\epsilon$-linearized differential}
\begin{equation}
\pa^{\epsilon}_{1}\colon\Aa_1(Y,\Lambda)/\Aa_{2}(Y,\Lambda)\to\Aa_1(Y,\Lambda)/\Aa_{2}(Y,\Lambda).
\end{equation}
\begin{dfn}
The {\em $\epsilon$-linearized contact homology} is the $\Z_{2}$-vector space
\begin{equation}
LCH_{\ast}(Y,\Lambda;\epsilon)=\krn(\pa^{\epsilon}_1)/\img(\pa^{\epsilon}_1).
\end{equation}
\end{dfn}
For simpler notation below we write
\[
Q(Y,\Lambda)=\Aa_1(Y,\Lambda)/\Aa_{2}(Y,\Lambda)
\]
and think of it as the graded vector space generated by the Reeb chords of $\Lambda$. Furthermore, the augmentation will often be clear from the context and we will drop it from the notation and write the differential as
\[
\pa_1\colon Q(Y,\Lambda)\to Q(Y,\Lambda).
\]

Consider an exact cobordism $(X,L)$ with $(+\infty)$-boundary $(Y^{+},\Lambda^{+})$ and $(-\infty)$-boundary $(Y^{-},\Lambda^{-})$.  Define the algebra map $\Phi\colon\Aa(Y^{+},\Lambda^{+})\to\Aa(Y^{-},\Lambda^{-})$ by mapping generators $c$ of $\Aa(Y^{+},\Lambda^{+})$ as follows
\[
\Phi(c)=\sum_{\dim(\M(c;\overline{b}))=0}|\M(c;\overline{b})|\,\overline{b},
\]
where $\overline{b}=b_1\dots b_k$ is a word of Reeb chord of $\Lambda^{-}$, where $\M(c;\overline{b})$ denotes the moduli space of holomorphic disks in $(X,L)$ with boundary on $L$, with positive puncture at $a$ and negative punctures at $b_1\dots b_k$. An argument completely analogous to the argument above showing $\pa^{2}=0$, looking at the boundary of $1$-dimensional moduli spaces shows that $\Phi\circ\pa^{+}=\pa^{-}\circ\Phi$, where $\pa^{\pm}$ is the differential on $\Aa(Y^{\pm},\Lambda^{\pm})$, i.e., that $\Phi$ is a chain map. Consequently, if $\epsilon^{-}\colon\Aa(Y^{-},\Lambda^{-})\to\Z_2$ is an augmentation then so is $\epsilon^{+}=\epsilon^{-}\circ \Phi$. In particular if $\Lambda^{-}=\emptyset$ and $(Y,\Lambda)=(Y^{+},\Lambda^{+})$ then $\Aa(Y^{-},\Lambda^{-})=\Z_{2}$ with the trivial differential and $\epsilon=\epsilon^{+}=\Phi$ is an augmentation of $\Aa(Y,\Lambda)$.

\subsection{Proof of Theorem \ref{t:iso}}\label{s:pft:iso}
If $Q_{[\alpha]}(Y,\Lambda)$ denotes the subspace of $Q(Y,\Lambda)$ generated by Reeb chords $c$ of action $\action(c)<\alpha$ then $Q_{[\alpha]}(Y,\Lambda)$ is a subcomplex of $Q(Y,\Lambda)$. By definition, the
map which takes a Reeb chord $c$ viewed as a generator of $\V_{[\alpha]}(X,L)$ to the dual $c^{\ast}$ of $c$ in the co-chain complex $Q_{[\alpha]}'(Y,\Lambda)$ of $Q_{[\alpha]}(Y,\Lambda)$ is an isomorphism intertwining the respective differentials. To prove the theorem it thus remains only to show that the monotonicity condition implies $H^{\ast}(Q'_{[\alpha]}(Y,\Lambda))=H^{\ast}(Q'(Y,\Lambda))$ for $\alpha>0$ large enough. This is straightforward: if $|c|=C_1\action(c)+C_0$ then
\[
H^{r}(Q'_{[\alpha]}(Y,\Lambda))=H^{r}(Q'(Y,\Lambda)),
\]
for $\alpha>\frac{r+1-C_0}{C_1}$.\qed

\section{Lagrangian Floer homology of exact cobordisms}
In this section we introduce a Lagrangian Floer cohomology of exact cobordisms. It is a generalization of the two copy version of the relative SFT of an exact cobordism $(X,L)$, $L=L_0\cup L_2$ to the case when each piece of $L$ is embedded but $L_0\cap L_1\ne \emptyset$. To prove that this theory has desired properties we will use a mixture of results from Floer homology of compact Lagrangian submanifolds and the SFT framework explained in Section \ref{S:relSFT}. After having set up the theory we state a conjectural lemma about how moduli spaces of holomorphic disks with boundary on $L\cup L'$, where $L'$ is a small perturbation of $L$, can be described in terms of holomorphic disks with boundary on $L$ and a version of Morse theory on $L$. We then show how Conjecture \ref{cnj:Seidel} and Corollary \ref{c:bigdiagram} follow from this conjectural description.

\subsection{The chain complex}
Let $X$ be a simply connected exact symplectic cobordism with $c_1(TX)=0$ and with good ends. Let $L_0$ and $L_1$ be exact Lagrangian cobordisms in $X$ with empty negative ends and with trivial Maslov classes. In other words, $(X,L_0)$ and $(X,L_1)$ are exact cobordisms with empty $(-\infty)$-boundaries. Let the $(+\infty)$-boundaries of $(X,L_{0})$ and $(X,L_{1})$ be $(Y,\Lambda_0)$ and $(Y,\Lambda_1)$, respectively.

Define
\[
C(X;L_0,L_1)=C_\infty(X;L_0,L_1)\oplus C_0(X;L_0,L_1)
\]
as follows. The summand $C_\infty(X;L_0,L_1)$ is the $\Z_2$-vector space of formal sums of Reeb chords which start on $\Lambda_0$ and end on $\Lambda_{1}$. The summand $C_0(X;L_0,L_1)$ is the $\Z_2$-vector space generated by the transverse intersection points in $L_0\cap L_1$.

In order to define grading and a differential on $C(X;L_0,L_1)$ we will consider the following three types of moduli spaces. The first type was considered already in \eqref{e:gradingCH}, if $a$ is a Reeb chord and $\overline{b}=b_1\dots b_k$ is a word of Reeb chords we write
\[
\M(a;\overline{b})
\]
for the moduli space of holomorphic curves in the symplectization $Y\times\R$ with boundary on $\Lambda_0\times\R\cup\Lambda_1\times\R$, with positive puncture at $a$ and negative punctures at $b_1,\dots,b_k$.

The second kind is the standard moduli spaces for Lagrangian Floer homology, if $x$ and $y$ are intersection points of $L_0$ and $L_1$ we write
\[
\M(x;y)
\]
for the moduli space of holomorphic disks in $X$ with two boundary punctures at which the disks are asymptotic to $x$ and $y$, with boundary on $L_0\cup L_1$, and which are such that in the orientation on the boundary induced by the complex orientation the incoming boundary component at $x$ maps to $L_0$.

Finally, the third kind is a mixture of these, if $c$ is a Reeb chord connecting $\Lambda_0$ to $\Lambda_1$ and if $y$ is an intersection point of $L_0$ and $L_1$ we write $\M(c;y)$ for the moduli space of holomorphic disks in $X$ with two boundary punctures at one the disk has a positive puncture at $c$ and at the other the disks is asymptotic to $y$.

If $a$ and $c$ are both Reeb chord generators then we define the grading difference between $a$ and $c$ to equal the formal dimension $\dim(\M(a;c))$. If $g$ is a Reeb chord- or an intersection point generator and if $x$ is an intersection point generator then we take the grading difference between $g$ and $x$ to equal $\dim(\M(a;c))+1$. Our assumptions on the exact cobordism guarantee that this is well defined.

\begin{figure}
\labellist
\small\hair 2pt
\pinlabel $a$   at 173 325
\pinlabel $c$   at 173 110
\pinlabel $b_1$ at 18 110
\pinlabel $b_2$ at 91 110
\pinlabel $e_1$ at 254 110
\pinlabel $e_2$ at 326 110
\pinlabel $e_3$ at 398 110
\pinlabel {\large $1$} at 170 186
\pinlabel {\large $0$} at 17 32
\pinlabel {\large $0$} at 90 32
\pinlabel {\large $0$} at 252 32
\pinlabel {\large $0$} at 324 32
\pinlabel {\large $0$} at 396 32
\endlabellist
\centering
\includegraphics[width=.5\linewidth]{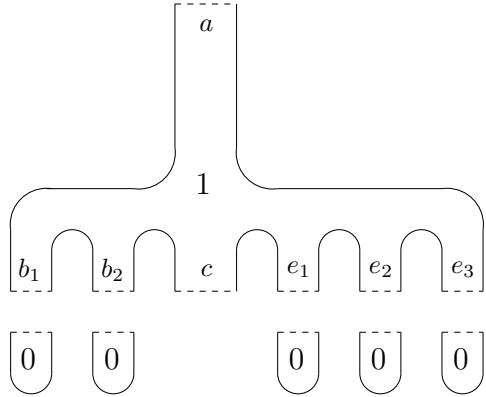}
\caption{A disk configuration contributing to the differential of $c$. Reeb chords are dashed. Numbers inside disk components indicate their dimension.}
\label{fig:cappeddisk}
\end{figure}

\begin{figure}
\labellist
\small\hair 2pt
\pinlabel $c$  at 50 221
\pinlabel $x$  at 50 29
\pinlabel $x$  at 322 25
\pinlabel $y$  at 322 111
\pinlabel {\large $0$} at 322 73
\pinlabel {\large $0$} at 50 110
\endlabellist
\centering
\includegraphics[width=.5\linewidth]{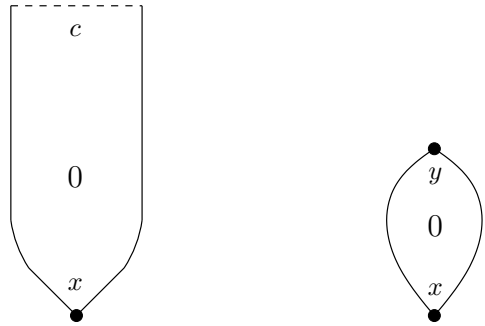}
\caption{Disks contributing to the differential of $x$. Reeb chords are dashed and double points appear as dots. Numbers inside disk components indicate their dimension.}
\label{fig:doubleptdisk}
\end{figure}

Write $C=C(X;L_0,L_1)$, $C_0=C_0(X;L_0,L_1)$, and $C_{\infty}=C_{\infty}(X;L_0,L_1)$. Define the differential $d\colon C\to C$, using the decomposition $C=C_\infty\oplus C_0$ to be given by the matrix
\[
d=
\left(
\begin{matrix}
d_{\infty} &  \rho\\
0     &   d_{0}
\end{matrix}
\right),
\]
where $d_{\infty}\colon C_\infty\to C_\infty$, $\rho\colon C_0\to C_\infty$, and $d_0\colon C_0\to C_0$ are defined as follows. Let $c$ be a Reeb chord from $\Lambda_0$ to $\Lambda_1$ and let $\theta\colon\Aa(Y,\Lambda_0)\to\Z_2$ and $\epsilon\colon\Aa(Y,\Lambda_1)\to\Z_2$ denote the augmentations induced by $L_0$ and $L_1$, respectively. Define
\begin{equation}\label{e:dinfty}
d_{\infty}(c)=\sum_{\dim(\M(a;\overline{b}\,c\,\overline{e}))=1}
|\widehat{\M}(a;\overline{b}\,c\,\overline{e})|
\epsilon(\overline{b})\theta(\overline{e})\, a,
\end{equation}
where $\overline{b}$ and $\overline{e}$ are words of Reeb chords from $\Lambda_0$ to $\Lambda_0$ and from $\Lambda_1$ to $\Lambda_1$, respectively, see Figure \ref{fig:cappeddisk}.
Let $x$ be an intersection point of $L_0$ and $L_1$. Define
\begin{equation}\label{e:rho}
\rho(x)=\sum_{\dim(\M(c;x))=0}
|\M(c;x)|\, c,
\end{equation}
where $c$ is a Reeb chord from $\Lambda_0$ to $\Lambda_1$, and
\begin{equation}\label{e:d0}
d_0(x)=\sum_{\dim(\M(y;x))=0}
|\M(y;x)|\, y,
\end{equation}
where $y$ is an intersection point of $L_0$ to $L_1$. See Figure \ref{fig:doubleptdisk}.
Then $d$ increases grading by $1$.

\begin{lma}
The map $d$ is a differential, i.e., $d^{2}=0$.
\end{lma}

\begin{pf}
We first check that $d_{\infty}^{2}=0$. Let $c$ and $a$ be a Reeb chords of index difference $2$. Consider two holomorphic disks in $\M(a';\overline{b}_-\,c\,\overline{e}_-)$ and $\M(a;\overline{b}_+\,a'\,\overline{e}_+)$ contributing to the coefficient of $a$ in $d^{2}(c)$. Gluing these two $1$-dimensional families at $a'$ and completing with Reeb chord strips at chords in $\overline{b}_+$ and $\overline{e}_+$, we find that the broken disk corresponds to one endpoint of a reduced moduli space $\widehat{\M}(a;\overline{b}\,c\,\overline{e})$, where $\overline{b}=\overline{b}_{+}\overline{b}_-$ and $\overline{e}=\overline{e}_-\overline{e}_{+}$. Note that there are three possible types of breaking at the boundary of $\widehat{\M}(a;\overline{b}\,c\,\overline{e})$:
\begin{itemize}
\item[(a)] breaking at a Reeb chord from $\Lambda_0$ to $\Lambda_1$,
\item[(b)] breaking at a Reeb chord from $\Lambda_0$ to $\Lambda_0$, and
\item[(c)] breaking at a Reeb chord from $\Lambda_1$ to $\Lambda_1$.
\end{itemize}
Summing the formal disks of all these boundary configurations give $0$ since the ends of a compact $1$-manifold cancels in pairs. To interpret this algebraically we define
\[
\tilde d_{\infty}\colon \Aa(Y,\Lambda_0)\otimes C_\infty\otimes\Aa(Y,\Lambda_1)\to
\Aa(Y,\Lambda_0)\otimes C_\infty\otimes\Aa(Y,\Lambda_1)
\]
as follows on generators:
\begin{equation*}
\tilde d_{\infty}(w_0\otimes c\otimes w_1)=\sum_{\dim(\M(a;\overline{b}\,c\,\overline{e}))=1}
|\widehat{\M}(a;\overline{b}\,c\,\overline{e})|
w_0\overline{b}\,a\,\overline{e}w_1.
\end{equation*}
Then the cancelation mentioned above implies with $\hat c=1\otimes c\otimes 1$ that
\[
\tilde d_{\infty}^{2}(\hat c)+\pa\otimes 1\otimes 1\,(\tilde d_{\infty}(\hat c))
+1\otimes 1\otimes \pa\,(\tilde d_{\infty}(\hat c))=0,
\]
where $\pa\colon\Aa(Y,\Lambda_j)\to\Aa(Y,\Lambda_j)$ denotes the contact homology differential. Here the first term corresponds to breaking of type (a), the second to type (c), and the third to type (b), respectively. By definition
\[
d_{\infty}(c)=\epsilon\otimes 1\otimes\theta\,(\tilde d_{\infty}(\hat c)).
\]
Consequently,
\begin{align*}
d_{\infty}^{2}(c)&=
\epsilon\otimes 1\otimes\theta\,(\tilde d_{\infty}^{2}(\hat c))\\
&=\epsilon\otimes 1\otimes\theta\,(\pa\otimes 1\otimes 1\,(\tilde d_{\infty}(\hat c))
+1\otimes 1\otimes \pa\,(\tilde d_{\infty}(\hat c)))=0,
\end{align*}
since $\epsilon\circ\pa=0$ and $\theta\circ \pa=0$.

The equation $d_0^{2}=0$ follows as in usual Lagrangian Floer homology from the fact that the terms contributing to $d_0^{2}$ are in $1-1$ correspondence with the ends of the $1$-dimensional moduli spaces of the form $\M(x;z)$, where $x$ and $z$ are intersection points.

Finally, to see that $d_\infty\circ \rho+\rho\circ d_0=0$ we consider $1$-dimensional moduli spaces of the form $\M(c,x)$ where $c$ is a Reeb chord from $\Lambda_0$ to $\Lambda_1$ and where $x$ is an intersection point. The analysis of breaking, see \cite{BEHWZ}, in combination with standard arguments from Lagrangian Floer theory shows that there are two possible breakings in the boundary of $\M(c,x)$:
\begin{itemize}
\item[(a)] breaking at an intersection point and
\item[(b)] breaking at Reeb chords.
\end{itemize}
In case (a) the broken configuration contributes to $\rho\circ d_0$. In case (b) the disk has two levels: the top level is a curve in a moduli space $\M(a;\overline{b}\,c\,\overline{e})$. The second level is a collection of rigid disks in $X$ with boundary on $L_0$ and $L_1$ and with positive punctures at Reeb chords in $\overline{b}$ and $\overline{e}$, respectively, and a rigid disk in $\M(c;x)$. Such a configuration contributes to $d_{\infty}\circ\rho$ and the desired equation follows.
\end{pf}

As above we let $\action(c)$ denote the action of a Reeb chord $c$. Define the action $\action(x)=0$ for intersection points $x\in L_0\cap L_1$. Then the differential on $C=C(X;L_0,L_1)$ increases action. Define $C_{[\alpha+]}\subset C$ as the subcomplex of formal sums in which all summands have action at least $\alpha$ and let $C_{[\alpha]}=C/C_{[\alpha+]}$ denote the corresponding quotient complex. Let $FH_{[\alpha]}^{\ast}(X;L_0,L_1)$ denote the homology of $C_{[\alpha]}$ and note that the natural projections give an inverse system of chain maps
\[
\pi^{\alpha}_{\beta}\colon C_{[\alpha]}\to C_{[\beta]},\quad \alpha>\beta.
\]
Define Lagrangian Floer cohomology $FH^{\ast}(X;L_0,L_1)$ as the inverse limit of the corresponding inverse system of homologies
\[
FH^{\ast}(X;L_0,L_1)=\underleftarrow{\lim}_{\alpha}\,FH^{\ast}_{[\alpha]}(X;L_0,L_1).
\]

\subsection{Chain maps and invariance}
Our proof of invariance of Lagrangian Floer cohomology $FH^{\ast}(X;L_0,L_1)$ under isotopies of $L_1$ uses three ingredients, homology isomorphisms induced by compactly supported isotopies, chain maps induced by joining cobordisms, and chain homotopies induced by compactly supported deformations of adjoined cobordisms. Before proving invariance, we consider these three separately.

\subsubsection{Compactly supported deformations}\label{ss:compactdef}
Let $(X,L_0)$ and $(X,L_1)$ be exact cobordisms as above. We first consider deformations of $(X,L_1)$ which are fixed in the positive end. More precisely, let $L_{1}^{t}$, $0\le t\le 1$ be a $1$-parameter family of exact Lagrangian cobordisms such that the positive end $\Lambda_{1}^{t}$ is fixed at $\Lambda_1$ for $0\le t\le 1$. Our proof of invariance is a generalization of a standard argument in Floer theory.

We first consider changes of the chain complex. Assuming that $L_{1}^{t}$ is generic there are a finite number of birth/death instances $0<t_1<\dots<t_m<1$ when two double points cancel or are born at a standard Lagrangian tangency moment. (At such a moment $t_j$ there is exactly one non-transverse intersection point $x\in L_0\cap L_1^{t_j}$, $\dim(T_xL_0\cap T_xL_1^{t_j})=1$, and if $v$ denotes the deformation vector field of $L_1$ at $x$ then $v$ is not symplectically orthogonal to $T_xL_0\cap T_xL_1^{t_j}$.)

For $0\le t\le 1$, let $\M^{t}$ denote a moduli space of the form $\M^{t}(g;x)$ or $\M^{t}(c)$, where $x$ is an intersection point, $c$ a Reeb chord, and $g$ either a Reeb chord or an intersection point, of holomorphic disks as considered in the definition of the differential, with boundary on $L_0$ and $L_1^{t}$. If $L^{t}_1$ is chosen generically then such a moduli space $\M^{t}$ is empty for all $t$ provided $\dim(\M^{t})<-1$ and there are a finite number of instances $0<\tau_1<\dots<\tau_k<1$ where there is exactly one disk of formal dimension $-1$ which is transversely cut out as a $0$-dimensional parameterized moduli space, see \cite[Lemma B.8]{E-JEMS}. We call these instances {\em $(-1)$-disk instances}. Furthermore, for generic $L_{1}^{t}$, birth/death instances and $(-1)$-disk instances are distinct.

For $I\subset[0,1]$, consider the parameterized moduli space of disks
\[
\M_{I}=\cup_{t\in I}\M^{t},
\]
where $\dim(\M^{t})=0$. If $I$ contains neither $(-1)$-disk instances nor birth/death instances then $\M_{I}$ is a $1$-manifold with boundary which consists of rigid disks in $\M^{t}$ and $\M^{t'}$ where $\pa I=\{t,t'\}$ and it follows by the definition of the differential that the chain complexes $C(X,L_0,L^{t}_{1})$ and $C(X;L_0,L^{t'}_{1})$ are canonically isomorphic. If on the other hand $I$ does contain $(-1)$-disk instances or birth/death instances then $\M_{I}$ has additional boundary points corresponding to broken disks at this instance. It is clear that in order to show invariance it is enough to show that the homology is unchanged over intervals containing only one critical moment. For simpler notation we take $I=[-1,1]$ and assume that there is a $(-1)$-disk instance or a birth/death at $t=0$.

We start with the case of a $(-1)$-disk. There are two cases to consider: either the $(-1)$-disk is mixed (i.e., has boundary components mapping both to $L_0$ and $L_1^{0}$) or it is pure (i.e., all of its boundary maps to $L_{0}^{0}$). Consider first the case of a mixed $(-1)$-disk which lies in a moduli space $\M(g;x)$, where $g$ is a Reeb chord or an intersection point and where $x$ is an intersection point. Write $C(-)=C(X;L_0,L^{-1}_{1})$ and $C(+)=C(X;L_0,L^{1}_{1})$, and let $d^{-}$ and $d^{+}$ denote the corresponding differentials.
\begin{lma}\label{l:chainisomixed}
Let $\phi\colon C(-)\to C(+)$ be the linear map defined on generators as follows
\[
\phi(g)=
\begin{cases}
x+y &\text{if }g=x,\\
g &\text{otherwise.}
\end{cases}
\]
If $\Phi(g)=\phi(g+d^{-}g)$ then $\Phi\colon C(-)\to C(+)$ is a chain isomorphism.
\end{lma}

\begin{pf}
We first note that the differentials $d^{+}$ and $d^{-}$ agree on the canonically isomorphic subspaces $C_{\infty}(-)$ and $C_{\infty}(+)$. In order to study the remaining part of the differential we consider the parameterized $1$-dimensional moduli space. Considering all punctures at intersection points as mixed, all disks which contribute to the differential are admissible, see Section \ref{s:formdi}, and since the $(-1)$-disk is mixed it follows from \cite[Lemma B.9]{E-JEMS} (or from a standard result in Floer theory, see \cite{Fl}, in case both generators are double points) that for any generator $g$ the following holds:
\[
d^{+} g =
\begin{cases}
d^- g +\phi(d^- g), &\text{if }|g|=|x|+1,\\
d^- g + x^{\ast}(g)\, d^{-}y, &\text{if }|g|=|x|,\text{ and}\\
d^{-} g, &\text{otherwise},
\end{cases}
\]
where $x^{\ast}$ is the dual covector of $x$. The lemma follows.
\end{pf}

Consider next the case of a pure disk. Here the situation is less straightforward. In order to get control of what is happening we need to introduce abstract perturbations of the $\bar\pa_J$-operator near the moduli space of $J$-holomorphic disks. We will give a short description here and refer to \cite[Section B.6]{E-JEMS} for details. As above, we study parameterized moduli spaces. In order to control contributions from broken disks of dimension $0$ with many copies of the $(-1)$-disk, a perturbation which time orders the complex structures at the negative pure $\Lambda_1$-punctures of any disk in the positive end is introduced. This perturbation needs to be extended over the entire moduli space of $1$-punctured holomorphic disks (below a fixed $(+)$-action) in the symplectization $(Y\times\R,\Lambda\times\R)$. Such a perturbation is defined energy level by energy level starting from the lowest one. The time ordering of the negative punctures implies that only one $(-1)$-disk at the time can be attached to any disk in the symplectization. Here it should be noted that the time ordering itself may introduce new $(-1)$-disks but that their positive punctures all lie close at $t=0$. More precisely, the count of $(-1)$-disks with positive puncture at a Reeb chord $c$ depends on the perturbation used on energy levels below $\action(c)$ and their positive puncture lie close to $t=0$ compared to the size of the time ordering perturbation of negative punctures mapping to $c$.

Let $\epsilon^-$ and $\epsilon^+$ denote the augmentations on $\Aa(Y,\Lambda)$ induced by $(X,L_{1}^{-1})$ and $(X,L_{1}^{1})$, respectively. It is a consequence of \cite[Lemma B.15]{E-JEMS} (which uses the perturbation scheme above) that there is a map $K\colon \Aa(Y,\Lambda_{1})\to\Z_2$ such that $K(c)$ counts $(-1)$-disks with positive puncture at $c$ and such that
\begin{equation}\label{e:augch}
\epsilon^-(c)+\epsilon^+(c)=\Omega_{K}(\pa c),
\end{equation}
here $\pa\colon\Aa(Y,\Lambda_{1})\to\Aa(Y,\Lambda_1)$ is the contact homology differential and if $w=b_1\dots b_m$ is a word of Reeb chords then
\[
\Omega_K(w)=\sum_j \epsilon^{-}(b_1)\dots\epsilon^{-}(b_{j-1})K(b_j)\epsilon^{+}(b_{j+1})\dots\epsilon^{+}(b_m).
\]

Define the linear map $\phi\colon C(-)\to C(+)$ as
\[
\phi(c)=\sum_{\dim(\M(a;\overline{b}\,c\,\overline{e}))=1}
|\widehat{\M}(a;\overline{b}\,c\,\overline{e})|
\Omega_{K}(\overline{b})\theta(\overline{e})\, a
\]
for generators $c\in C_{\infty}(-)$ and $\phi(x)=0$ for generators $x\in C_{0}(-)$.

\begin{lma}\label{l:chainisopure}
The map $\Phi\colon C_-\to C_+$
\[
\Phi(c)=c+\phi(c),
\]
is a chain isomorphism.
\end{lma}

\begin{pf}
The map is an isomorphism since the action of any chord in $\phi(c)$ is larger than that of $c$. We thus need only show it is a chain map or in other words that
\[
d^+ + d^{-}=d^{+}\circ\phi + d^{-}\circ \phi.
\]

Consider first the operator in the left hand side acting on a Reeb chord $c$. According to \eqref{e:augch}, broken disk configurations which contributes to $d^+(c)+d^-(c)$ are of the following form, see Figure \ref{fig:purehs0}:

\begin{figure}
\labellist
\small\hair 2pt
\pinlabel $a$  at 225 325
\pinlabel $c$  at 225 216
\pinlabel $b_s$  at 155 178
\pinlabel {\large $1$} at 225 251
\pinlabel {$-1$} at 155 25
\pinlabel {\large $1$} at 155 105
\pinlabel {\large $0$} at 10 105
\pinlabel {\large $0$} at 46 105
\pinlabel {\large $0$} at 82 105
\pinlabel {\large $0$} at 226 105
\pinlabel {\large $0$} at 262 105
\pinlabel {\large $0$} at 298 105
\pinlabel {\large $0$} at 334 105
\pinlabel $\epsilon^{-}$ at 118 25
\pinlabel $\epsilon^{-}$ at 83 25
\pinlabel $\epsilon^{-}$ at 47 25
\pinlabel $\epsilon^{-}$ at 12 25
\pinlabel $\epsilon^{+}$ at 191 25
\pinlabel $\theta$ at 261 25
\pinlabel $\theta$ at 299 25
\pinlabel $\theta$ at 333 25
\endlabellist
\centering
\includegraphics[width=.8\linewidth]{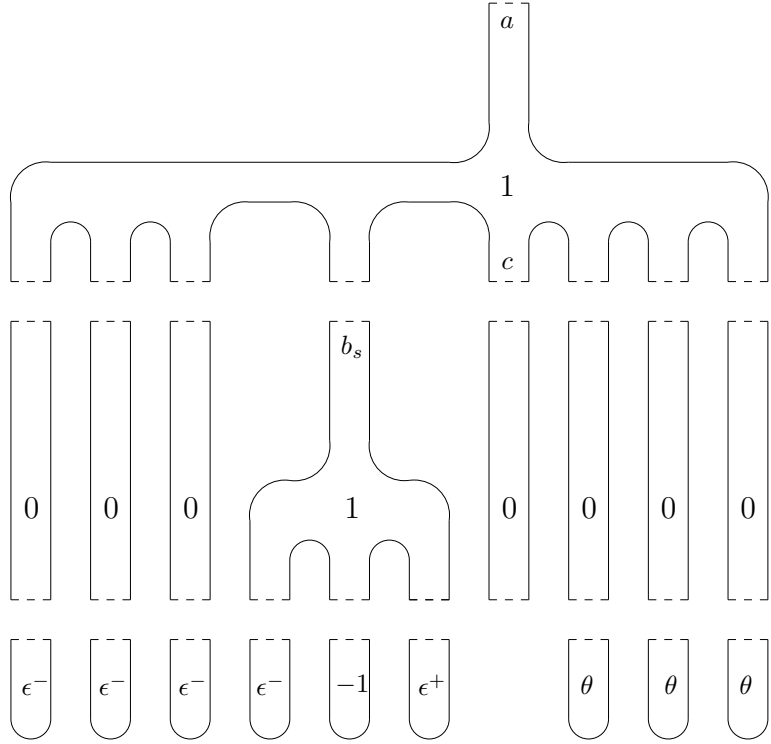}
\caption{A disk contributing to $d^{+}(c)-d^{-}(c)$: the disk breaks at the pure Reeb chord $b_s$ and the $(-1)$-disk is attached to one of the negative ends of the disk with positive puncture at $b_s$.}
\label{fig:purehs0}
\end{figure}

\begin{itemize}
\item[$(1)$] The top level is a $1$-dimensional disk in $Y\times\R$ with positive puncture at a Reeb chord $a$ connecting $\Lambda_0$ to $\Lambda_1$, followed by $k$ negative punctures at Reeb chords $b_1\dots b_k$ connecting $\Lambda_1$ to itself, followed by a negative puncture at $c$, in turn followed by $r$ negative punctures at Reeb chords $e_1,\dots e_r$ connecting $\Lambda_1$ to itself.
\item[$(2)$] The middle level consists of Reeb chord strips at all negative punctures except $b_s$ for some $1\le s\le k$. At $b_s$ a $1$-dimensional disk in $Y\times\R$ with boundary on $\Lambda_1\times\R$ and with negative punctures at $q_1,\dots,q_m$ is attached.
\item[$(3)$] The bottom level consists of rigid disks with boundary on $L_0$ and positive puncture at $e_j$ attached at all punctures $e_j$, $j=1,\dots,r$, rigid disks with boundary on $L_{1}^{-1}$ attached at punctures $b_j$, $1\le j\le s-1$ and at punctures $q_j$, $1\le j\le v$, a $(-1)$-disk attached at $q_v$, and rigid disks with boundary on $L_{1}^{1}$ at punctures $b_j$, $s+1\le j\le k$.
\end{itemize}
Consider gluing of the top and middle level above in the symplectization. This gives one boundary component of a reduced $1$-dimensional moduli space. The other boundary component corresponds to one of three breakings:  at a pure $\Lambda_1$-chord, at a chord connecting $\Lambda_0$ to $\Lambda_1$, or breaking at a pure $\Lambda_0$ chord. The first type of breaking contributes to $d_{-}(c)+d_{+}(c)$ as well, see Figure \ref{fig:purehs0}. The second type contributes to either $d^{+}\circ\phi$ or $\phi\circ d^{-}$ depending on which factor the $(-1)$-end goes, see Figures \ref{fig:purehs1A} and \ref{fig:purehs1B}, respectively, and finally, the total contribution of the third type of breaking equals $0$ since  $\theta\circ \pa=0$, where $\pa\colon\Aa(Y,\Lambda_0)\to\Aa(Y,\Lambda_0)$ is the contact homology differential and where $\theta\colon \Aa(Y,\Lambda_0)\to\Z_2$ is the augmentation, see Figure \ref{fig:purehs2}.

\begin{figure}
\labellist
\small\hair 2pt
\pinlabel $a$  at 225 325
\pinlabel {\large $1$} at 225 251
\pinlabel {\large $0$} at 225 107
\pinlabel {\large $0$} at 297 107
\pinlabel {\large $0$} at 333 107
\pinlabel {\large $0$} at 81 107
\pinlabel {\large $0$} at 45 107
\pinlabel {\large $0$} at 9 107
\pinlabel $c$ at 226 182
\pinlabel {$-1$} at 153 25
\pinlabel $\epsilon^{-}$ at 118 25
\pinlabel $\epsilon^{-}$ at 83 25
\pinlabel $\epsilon^{-}$ at 47 25
\pinlabel $\epsilon^{-}$ at 12 25
\pinlabel $\epsilon^{+}$ at 191 25
\pinlabel $\theta$ at 261 25
\pinlabel $\theta$ at 299 25
\pinlabel $\theta$ at 333 25
\endlabellist
\centering
\includegraphics[width=.8\linewidth]{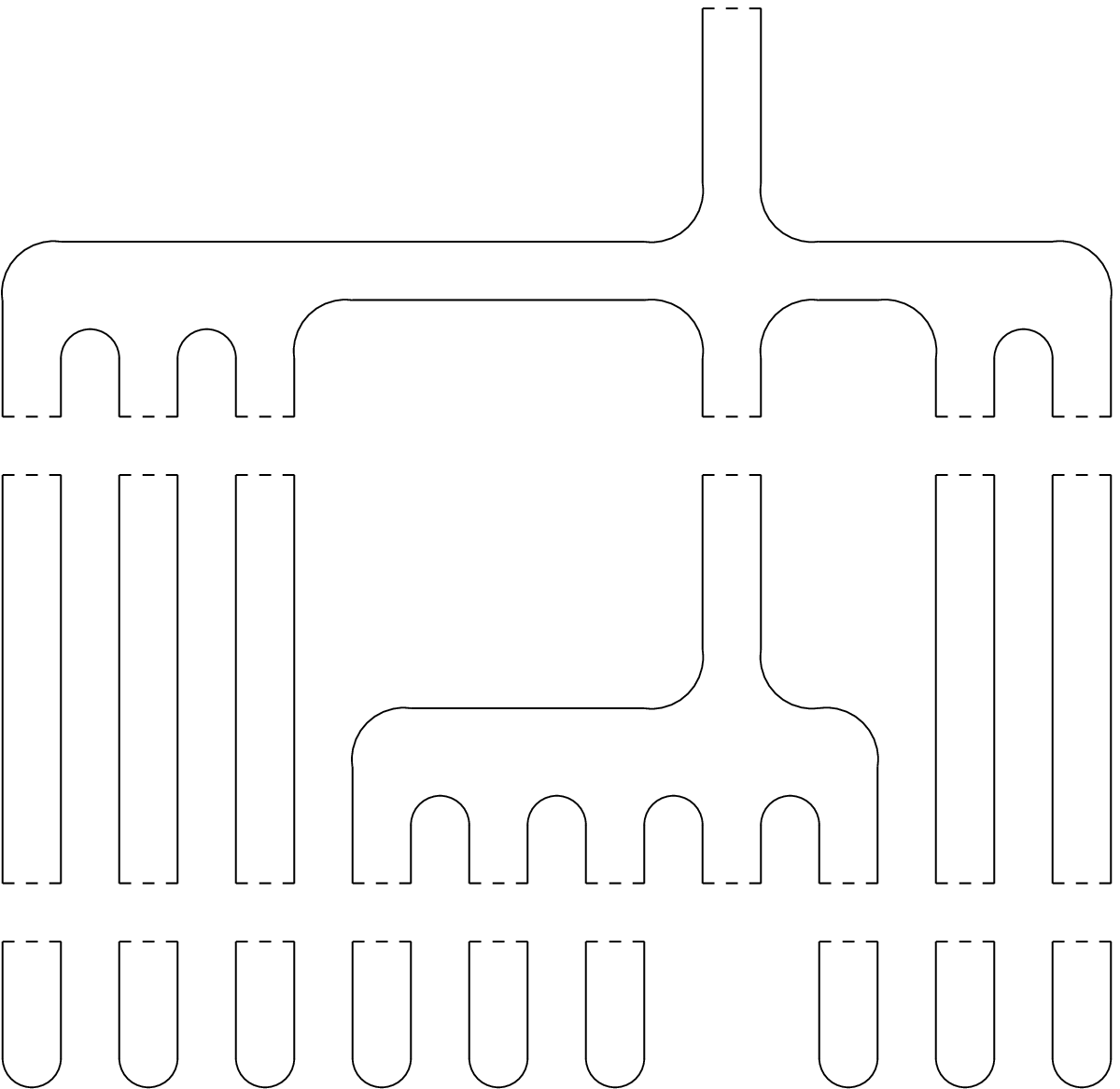}
\caption{A disk contributing to $d^{+}\circ\phi(c)$: the disk breaks at the mixed chord $c$ and the $(-1)$-disk is attached to the disk with positive puncture at $c$.}
\label{fig:purehs1A}
\end{figure}

\begin{figure}
\labellist
\small\hair 2pt
\pinlabel $a$  at 225 325
\pinlabel {\large $0$} at 225 251
\pinlabel {\large $1$} at 225 107
\pinlabel {\large $0$} at 297 107
\pinlabel {\large $0$} at 333 107
\pinlabel {\large $0$} at 81 107
\pinlabel {\large $0$} at 45 107
\pinlabel {\large $0$} at 9 107
\pinlabel $c$ at 226 182
\pinlabel $\epsilon^{+}$ at 153 25
\pinlabel $\epsilon^{+}$ at 118 25
\pinlabel $\epsilon^{+}$ at 83 25
\pinlabel {$-1$} at 46 25
\pinlabel $\epsilon^{-}$ at 12 25
\pinlabel $\epsilon^{+}$ at 191 25
\pinlabel $\theta$ at 261 25
\pinlabel $\theta$ at 299 25
\pinlabel $\theta$ at 333 25
\endlabellist
\centering
\includegraphics[width=.8\linewidth]{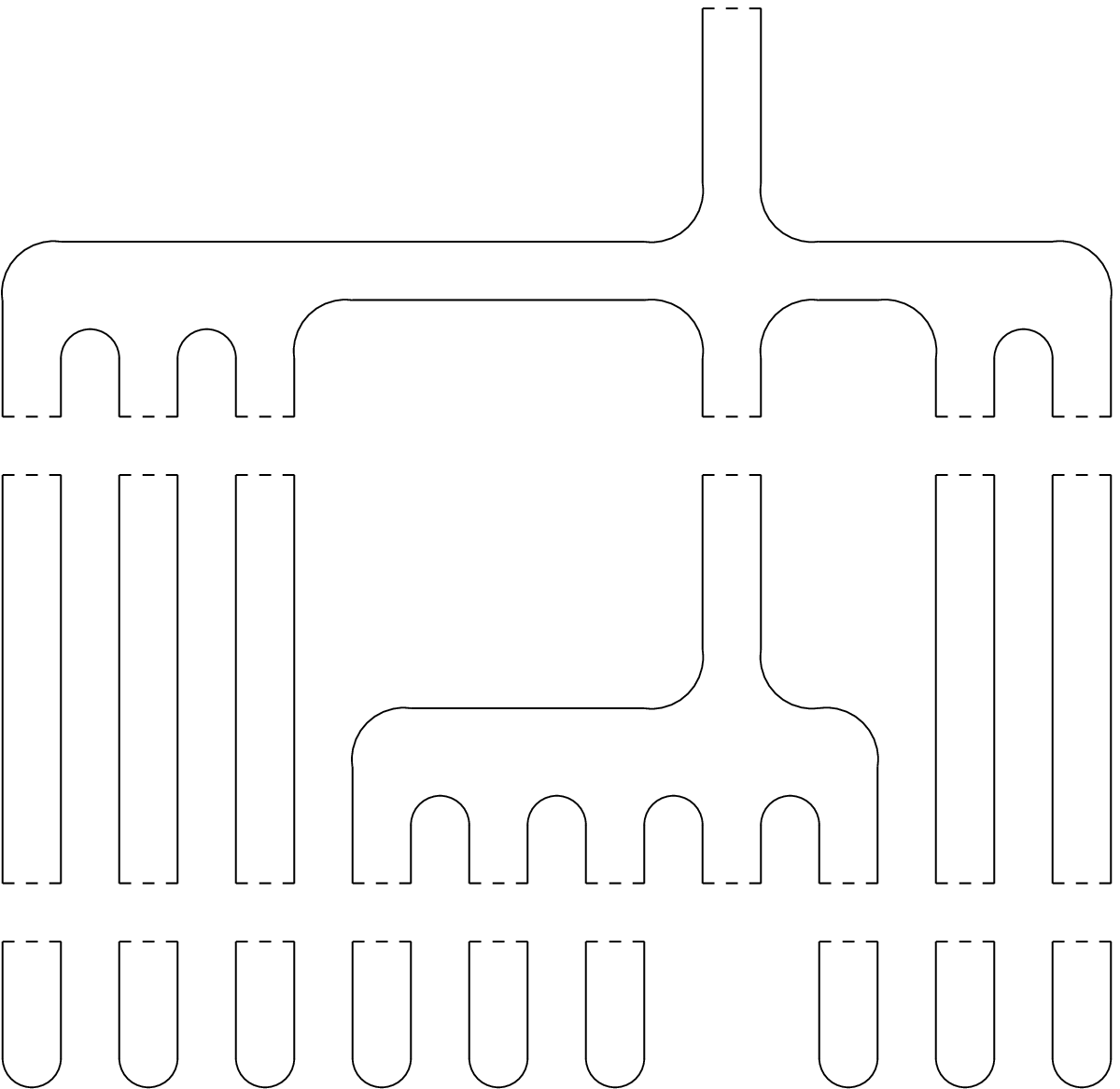}
\caption{A disk contributing to $\phi\circ d^{-}(c)$: the disk breaks at the mixed chord $c$ and the $(-1)$-disk is attached to the disk with positive puncture at $a$.}
\label{fig:purehs1B}
\end{figure}

\begin{figure}
\labellist
\small\hair 2pt
\pinlabel $a$  at 225 325
\pinlabel {\large $0$} at 225 251
\pinlabel {\large $0$} at 225 107
\pinlabel {\large $1$} at 297 107
\pinlabel {\large $0$} at 117 107
\pinlabel {\large $0$} at 154 107
\pinlabel {\large $0$} at 190 107
\pinlabel {\large $0$} at 81 107
\pinlabel {\large $0$} at 45 107
\pinlabel {\large $0$} at 9 107
\pinlabel $c$ at 226 182
\pinlabel $\epsilon^{+}$ at 153 25
\pinlabel $\epsilon^{+}$ at 118 25
\pinlabel {$-1$} at 83 25
\pinlabel $\epsilon^{-}$ at 46 25
\pinlabel $\epsilon^{-}$ at 12 25
\pinlabel $\epsilon^{+}$ at 191 25
\pinlabel $\theta$ at 261 25
\pinlabel $\theta$ at 299 25
\pinlabel $\theta$ at 333 25
\endlabellist
\centering
\includegraphics[width=.8\linewidth]{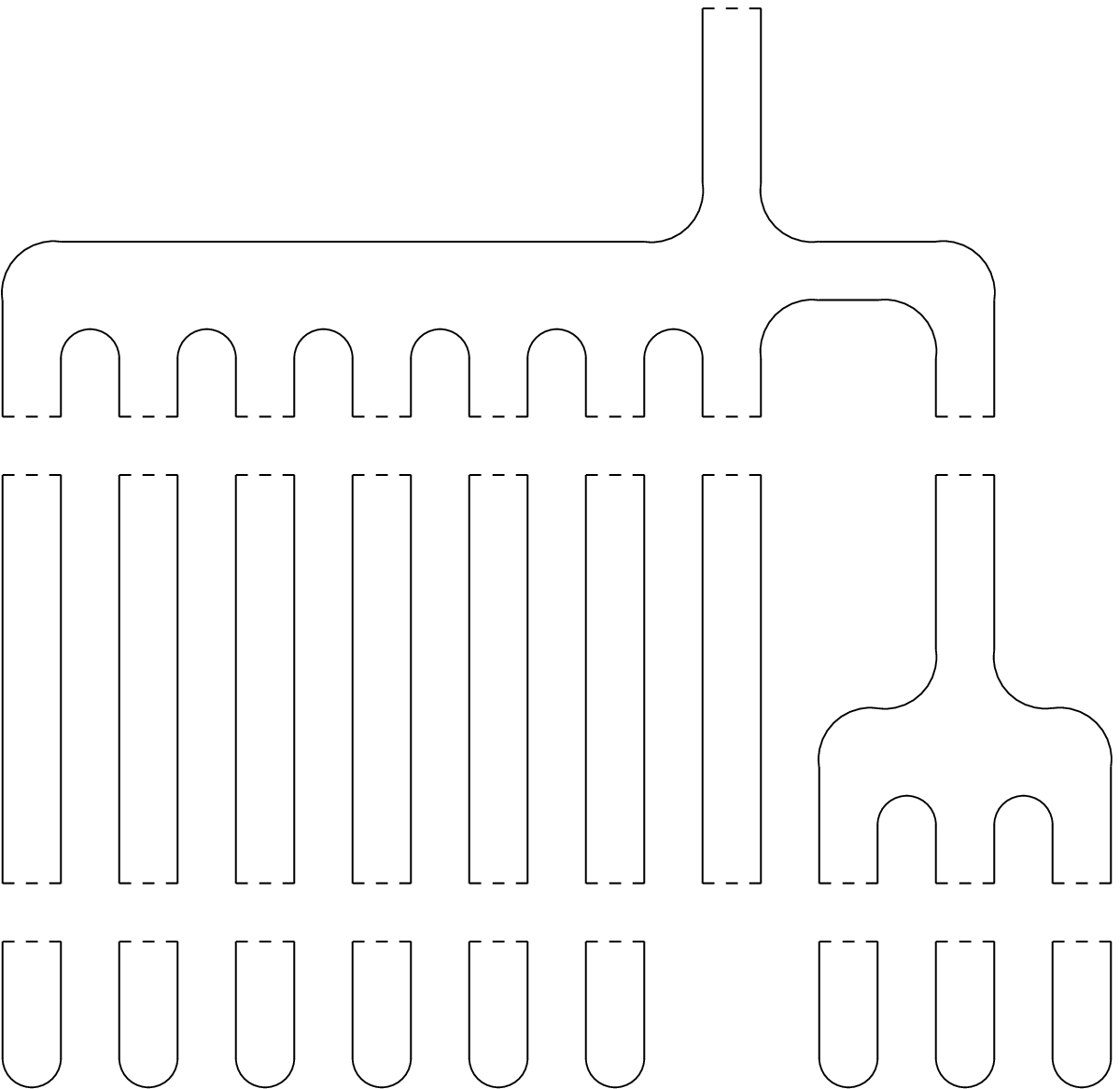}
\caption{Disks with total contribution $0$: the disk breaks at the pure $\Lambda_0$ chord $e$ and the total contribution vanishes since $\theta\circ\pa=0$.}
\label{fig:purehs2}
\end{figure}

Second consider the operators acting on a double point $x$. The chain map property in this case follows from an argument similar to the one just given. Additional boundary components of the parameterized moduli space $\M^{t}(c;x)$ corresponds to two level disks with top level a disk as in $(1)$ above and with bottom level as in $(3)$ above with the addition that there is a rigid disk with positive puncture at $c$ and a puncture at $x$. We conclude that $\Phi$ is a chain map.
\end{pf}

Finally, consider a birth/death moment involving intersection points $x$ and $y$. Then $dx=y+v$ where $v$ does not contain any $y$-term, see \cite{Fl}. Define the map $\Phi\colon C(+)\to C(-)$
by
\[
\Phi(x)=0,\quad \Phi(y)=v,\quad\text{and}\quad\Phi(w)=w \text{ for }w\ne x,y.
\]
Define the map $\Psi\colon C(-)\to C(+)$ by
\[
\Psi(c)=c+y^{\ast}(d^+ c)x.
\]

\begin{lma}\label{l:chainisob/d}
The maps $\Phi$ and $\Psi$ are chain maps which induce isomorphisms on homology.
\end{lma}
\begin{pf}
A straightforward generalization of the gluing theorem \cite[Proposition 2.16]{EES-JDGanalysis} (in the case when only one disk is glued at the degenerate intersection) shows that if $g$ is a generator of $C(-)$ then
\[
d^-(c)=d^+(c)+y^{\ast}(d^+(c))v.
\]
The lemma then follows from a straightforward calculation.
\end{pf}

\subsubsection{Joining cobordisms}\label{ss:joincob}
Let $(X,L_0)$ and $(X,L_1)$ be cobordisms as above. Consider the trivial cobordism $(Y\times\R,\Lambda_0\times\R)$ and some cobordism $(Y\times\R,L^{a}_1)$ where the $(-\infty)$-boundary of $L_1^{a}$ equals $\Lambda_{1}$. Then we can join these cobordisms to $(X,L_0)$ and $(X,L_1)$, respectively. This results in a new pair of cobordisms $(X,L_0)$ and $(X,\tilde L_1)$. We define a map
\[
\Phi\colon C(X;L_0,L_1)\to C(X;L_0,\tilde L_1,)
\]
as follows.

For points $x$ in $L_0\cap L_1$,
\[
\Phi(x)=x.
\]
For Reeb chords $c$ from $\Lambda_0$ to $\Lambda_1$,
\begin{align*}
\Phi(c)=&\sum_{\dim\M(x;\overline{b}\, c\,\overline{e})=0}
|\M(x;\overline{b}\,c\,\overline{e})|\epsilon(\overline{b})\theta(\overline{e})\, x\\
+&\sum_{\dim\M(a;\overline{b}\,c\,\overline{e})=0}
|\M(a;\overline{b}\,c\,\overline{e})|\epsilon(\overline{b})\theta(\overline{e})\, a
\end{align*}
where the first sum ranges over intersection points $x$ in $(\Lambda_0\times\R)\cap L^{a}_1$ and the second over Reeb chords $a$ connecting $\Lambda_0$ to $\Lambda_1^{a}$.

\begin{lma}\label{l:cobchmap}
The map $\Phi$ is a chain map. That is, if $d$ and $\tilde d$ denote the differentials on $C(X,L_0, L_1)$ and $C(X,L_0,\tilde L_1)$, respectively, then
\[
\tilde d\circ \Phi=\Phi\circ d.
\]
\end{lma}

\begin{pf}
Consider first an intersection point $x\in L_0\cap L_1$. A moduli space contributing to  $\tilde d(x)$ is of the form $\M(g;x)$, where $g$ is either at double point or a Reeb chord connecting $\Lambda_0$ to $\tilde \Lambda_1$. Consider stretching along the hypersurface where the cobordisms are joined. It is a consequence of \cite{BEHWZ} that families of disks in $\M(g;x)$ converge to broken disks in the limit. Since $\dim(\M(g;x))=0$ every level in the limit must have dimension $0$ by transversality. Now if the disk does not break then $g=y$ is another intersection point of $L_0\cap L_1$ and since $\Phi(y)=y$, it contributes to $\Phi\circ d$ as well. If on other hand the disk does break then $g=c$ or $g=z$, where $c$ is a Reeb chord connecting $\Lambda_{0}$ to $\tilde\Lambda_1$ and $z$ is an intersection point in $(\Lambda_0\times\R)\cap L^{a}_1$. In particular, by admissibility of the disks in $\M(g;x)$ there is exactly one Reeb chord $c'$ connecting $\Lambda_0$ to $\Lambda_1$ in the broken disk which hence contributes to $\Phi\circ d$ by definition of $\phi(c')$.

Consider second a Reeb chord $c$ connecting $\Lambda_0$ to $\Lambda_1$. To show that the chain map equation holds in this case we consider $1$-dimensional moduli spaces $\M(a;\overline{b}\,c\,\overline{e})$ where $\overline{b}=b_1\dots b_k$ and $\overline{e}=e_1\dots e_r$ are words of Reeb chords connecting $\Lambda_{1}$ to itself and connecting $\Lambda_0$ to itself, respectively, and such that $|b_j|=0$ all $j$ and $|e_l|=0$ all $l$. The boundary of such a moduli space consists of $2$-level disks: one level in $(Y\times\R,L_0\cup L^{a}_1)$ and one in either the positive- or the negative end of this symplectic manifold. If the $1$-dimensional disk lies in the positive end then configuration contributes to $\tilde d\circ \Phi$. If the $1$-dimensional disk lies in the negative end there are two cases to consider: either the positive puncture of this piece is a mixed Reeb chord or it is pure. In the former case the disk contributes to $\Phi\circ d$, the total contribution from disk of the later type vanishes since $\theta\circ\pa=0$ and $\epsilon\circ\pa=0$.
\end{pf}

\subsubsection{Joining maps and deformations}\label{ss:joinmapdef}
We next consider generic $1$-parameter families of cobordisms as considered in Section \ref{ss:joincob}. More precisely, let $(X,L_0)$ and $(X,L_1)$ be exact cobordisms as usual. Let $(Y\times\R,L^{a}_{1}(t))$ be a $1$-parameter family of exact cobordisms which is constant outside a compact set and such that the $(-\infty)$-boundary of $L^{a}_{1}(t)$ equals $\Lambda_1$. Adjoining $(Y\times\R,L^{a}_{1}(t))$ to $(X,L_1)$, we obtain a $1$-parameter family of cobordisms $(X,\tilde L_1(t))$ and corresponding chain maps
\[
\Phi_t\colon C(X;L_0,L_1)\to C(X;L_0,\tilde L_1(t)).
\]
Furthermore, since $(X,\tilde L_1(t))$ and $(X,\tilde L_1(t'))$ are related by a compact deformation Lemmas \ref{l:chainisomixed}, \ref{l:chainisopure}, and \ref{l:chainisob/d} provide chain maps
\[
\Psi_{tt'}\colon C(X;L_0,L_1(t))\to C(X;L_0,\tilde L_1(t'))
\]
which induce isomorphisms on homology. In fact, unless the interval between $t$ and $t'$ contain $(-1)$-disk instances or birth/death instances, the map $\Psi_{tt'}$ is the canonical identification map on generators. Here the birth/deaths take place in the added cobordism and the $(-1)$-disk instances correspond to $(-1)$-disk instances in the added cobordism. As we shall see below, if the interval between $t$ and $t'$ contains a birth/death- or a $(-1)$-disk instance then there is a chain homotopy connecting the chain maps $\Psi_{tt'}\circ\Phi_t$ and $\Phi_{t'}$. The proofs of these results are similar to the proofs of results in Sections \ref{ss:compactdef} and \ref{ss:joincob} and many details from there will not be repeated. For convenient notation below we take $t=-1$, $t'=1$ and assume that the critical instance is at $t=0$. Furthermore, we write $C=C(X,L_0,L_1)$ with differential $d$, $\tilde C_{\pm}=C(X,L_0,\tilde L_1({\pm1}))$ with differential $\tilde d_{\pm}$, $\Phi_{\pm}=\Phi_{\pm 1}$, and $\Psi=\Psi_{-1\, 1}$

Consider first the case of a pure $(-1)$-disk in $(Y\times\R,L^{a}_{1}(t))$. Applying our perturbation scheme which time orders the negative punctures of disks in the positive end $(Y\times\R,\Lambda_{1}^{a})$, such a disk give rise to several $(-1)$-disks in $(Y\times\R,L^{a}_{1}(t))$. The $(-1)$-disks in the total cobordism $(X,\tilde L_1(t))$ are then these $(-1)$-disks in $(Y\times\R,L^{a}_{1}(t))$, capped off with $0$-dimensional rigid disks in $(X,L_1)$, see \cite[Lemmas 4.3 and 4.4]{E-JEMS}.

\begin{lma}\label{l:chhmtpypure(-1)}
If there is a pure $(-1)$-disk at $t=0$, then the following diagram commutes
\[
\begin{CD}
C @>{\Phi_-}>> \tilde C_-\\
@V{\id}VV @VV{\Psi}V\\
C @>>{\Phi_+}> \tilde C_+
\end{CD}\quad\quad.
\]
\end{lma}

\begin{pf}
Consider the parameterized $1$-dimensional moduli space corresponding to a $0$-dimensional moduli space contributing to $\Phi_+$. A boundary component at $t=1$ contributes to $\Psi\circ\Phi_-$. A boundary component in the interior of $[-1,1]$ is a broken disk consisting of a $1$-dimensional disk in an end and a $(-1)$-disk and $0$-disks in the cobordisms. If the $1$-dimensional disk lies in the upper end then the broken disk contributes to $\Psi\circ\Phi_-$ and, as usual, the total contribution of disks with $1$-dimensional disk in the lower end equals $0$ since augmentations are chain maps. 
The result follows.
\end{pf}

Second, consider the case of a mixed $(-1)$-disk. By admissibility any disk contributing to the chain maps then contain at most one such $(-1)$-disk, see \cite[Lemma 2.8]{E-JEMS}. Define the map $K\colon C\to \tilde C_+$ as follows
\[
K(x)=0
\]
if $x$ is an intersection point generator and
\[
K(c)=\sum_{\dim(\M_{I}(g;\overline{b}\, c\,\overline{e}))=0}|\M_I(g;\overline{b}\, c\,\overline{e})|\epsilon(\overline{b})\theta(\overline{e})\,g,
\]
if $c$ is a Reeb chord generator, where $\M_{I}(g;\overline{b}\, c\,\overline{e})$ is the parameterized moduli space of disks in $Y\times\R$ with boundary on $\Lambda_0\times\R$ and $L_1^{a}(t)$, and where $\overline{b}$ and $\overline{e}$ are Reeb chords of $\Lambda_{0}$ and $\Lambda_{1}$, respectively.
\begin{lma}\label{l:chhmtpymixed(-1)}
If there is a mixed $(-1)$-disk at $t=0$, then chain maps in the following diagram
\[
\begin{CD}
C @>{\Phi_-}>> \tilde C_-\\
@V{\id}VV @VV{\Psi}V\\
C @>>{\Phi_+}> \tilde C_+
\end{CD}
\]
satisfies $\Psi\circ\Phi_-+\Phi_+= K\circ d + \tilde d_+\circ K$
\end{lma}

\begin{pf}
Consider first the left hand side acting on an intersection point $x\in L_0\cap L_1$. Both maps $\Phi_-$ and $\Phi_+$ are then inclusions and by definition, $\Psi(x)$ counts $(-1)$-disks in $\M(g;x)$ where $g$ is a generator. Since there are no $(-1)$-disks in the lower cobordism we find that in a moduli space which contributes to $\Psi(x)$ the generator $g$ is either an intersection point or a Reeb chord in the upper cobordism. Consider now the splittings of a $(-1)$-disk as we stretch over the joining hypersurface: it splits into a $(-1)$-disk in the upper cobordism and $0$-dimensional disks in the lower. By definition, the count of such split disks equal $K(d(x))$. Since $K(x)=0$, the chain map equation follows.

Consider second the left hand side acting on a Reeb chord generator $c$. Here $\Phi_\pm(c)$ are given by counts of $0$-dimensional disks in the upper cobordism and $\Psi$ is a count of $(-1)$-disks emanating at double points. In particular $\Psi(a)=a$ for any Reeb chord. Consider now a moduli space which contributes to $\Phi_+$. The boundary of the corresponding parameterized moduli space consists of rigid disks over endpoints as well as broken disks with a $1$-dimensional disk in either symplectization end and a $(-1)$-disk in the cobordism. The total count of such disks gives the desired equation after observing that splittings at pure chords do not contribute for the usual reason.
\end{pf}

Third, consider a birth/death instance.
\begin{lma}\label{l:chhmtpyb/d}
If there is a birth/death instance at $t=0$, then the following diagram commutes
\[
\begin{CD}
C @>{\Phi_-}>> \tilde C_-\\
@A{\id}AA @AA{\Psi}A\\
C @>>{\Phi_+}> \tilde C_+
\end{CD}\quad\quad.
\]
\end{lma}

\begin{pf}
Assume that the canceling pair of double points is $(x,y)$ with $\tilde d_{+}x=y+v$ as in Lemma \ref{l:chainisob/d}. As there, disks from $x$ to $v$ can on the one hand be glued to rigid disks ending at $y$ resulting in rigid disks, and on the other to rigid disks ending at $x$ resulting in non-rigid disks. Commutativity then follows from a straightforward calculation.
\end{pf}

\subsubsection{Invariance}
Let $(X,L_0)$ and $(X,L_1)$ be exact cobordisms as above.
\begin{thm}\label{t:invFH}
The Lagrangian Floer homology $FH^{\ast}(X;L_0,L_1)$ is invariant under exact deformations of $L_1$.
\end{thm}

\begin{pf}
The proof is similar to the proof of Theorem \ref{t:rSFTinv}:
any deformation considered can be subdivided into a compactly supported deformation and a Legendrian isotopy at infinity. The former type induces isomorphisms on homology by Lemmas \ref{l:chainisomixed}, \ref{l:chainisopure}, and \ref{l:chainisob/d}. The later type of deformation gives rise to an invertible exact cobordism which in turn gives a chain map on homology by Lemma \ref{l:cobchmap}. Lemmas \ref{l:chhmtpypure(-1)}, \ref{l:chhmtpymixed(-1)}, and \ref{l:chhmtpyb/d} show that on the homology level these maps are independent of compact deformations of the cobordism added. A word by word repetition of the proof of Theorem \ref{t:rSFTinv} then finishes the proof.
\end{pf}

\subsection{Holomorphic disks for $FH^{\ast}(X;L,L')$}
Let $(X,L)$ be an exact cobordism as considered above. Let $L'$ denote a copy of $L$. In order to make $L$ and $L'$ transverse, we identify a neighborhood of $L\subset X$ with the cotangent bundle $T^{\ast}L$. Pick a function $\tilde F\colon L\to \R$ such that
\[
\tilde F(y,t)=C+t,
\]
where $C$ is a constant, for $(y,t)\in \Lambda\times[T,\infty)$ for some $T>0$. Cut the function $F$ off outside a small neighborhood of $L$ and let $L''$ be the time $1$-flow of the Hamiltonian of $\epsilon F$ for small $\epsilon$. Then $L$ intersects $L''$ transversely. However, in the end where $F=t$, the Hamiltonian just shifts $\Lambda$ along the Reeb flow so that there is a Reeb chord from $\Lambda$ to $\Lambda''$ at every point of $\Lambda$. In order to perturb out from this Morse-Bott situation, we identify a neighborhood of $\Lambda\subset Y$  with $J^{1}(\Lambda)$, fix a Morse function $f\colon\Lambda\to\R$ and let $\Lambda'$ denote the graph of the $1$-jet extension of $f$. Then $\Lambda$ and $\Lambda'$ are contact isotopic via the contact isotopy generated by the time dependent contact Hamiltonian $H_{t}(q,p,z)=\psi(t)f(q)$, $(q,p,z)\in T^{\ast}\Lambda\times\R$ and where $\psi(t)$ is a cut off function. Take $f$ and $\frac{d\psi}{dt}$ very small and adjoin the cobordism $Y\times\R$ with the symplectic form $d e^{t}(\lambda-H_t)$ to $(X,L)$. This gives the desired $L'$ with $\Lambda'$ as $(+\infty)$-boundary.

We will state a conjectural lemma which gives a description of holomorphic disks with boundary on $L\cup L'$. To this end we first describe a version of Morse theory on $L$ and second discuss intersection points in $L\cap L'$ and Reeb chords of $\Lambda\cup \Lambda'$.

Consider the gradient equation $\dot x=\nabla F(x)$ of the function $F\colon L\to\R$, where
\[
F(x)=\tilde F(x)+\psi(t)f(y),
\]
where we write $x=(y,t)\in\Lambda\times[T,\infty)$. It is easy to see that if a solution of $\dot x=\nabla F(x)$ leaves every compact then it is exponentially asymptotic to a {\em critical point solution} of the form $t\mapsto (y_0,t)$, where $\nabla f(y_0)=0$ in $\Lambda\times[0,\infty)$. Furthermore, every sequence of solutions of $\dot x=\nabla F(x)$ has a subsequence which converges to a several level solution with one level in $L$ and levels in $\Lambda\times\R$ which are solutions to the gradient equation of $f+t$, asymptotic to critical point solutions. We call solutions which have formal dimension $0$ (after dividing out re-parametrization) and which are transversely cut out {\em rigid flow lines}. Consider a moduli space $\M(c)$ of holomorphic disks in $X$ with boundary on $L$ or a moduli space $\M(c;\overline{b})$ in $Y\times\R$ with boundary on $\Lambda\times\R$. Marking a point on the boundary there is an evaluation map $\ev\colon \M^{\ast}(c)\to L$ or $\ev\colon\M^{\ast}(c;\overline{b})\to\Lambda\times\R$, see \cite[Section 6]{EESa-duality}. Consider now a flow line of $F$ or of $f+t$ which hits the image of $\ev$. We call such a configuration a {\em generalized disk} and we say that it is rigid if it has formal dimension $0$ (after dividing out the $\R$-translation in $Y\times\R$) and if it is transversely cut out.

Note that points in $L\cap L'$ are in $1-1$ correspondence with critical points of $F$. We use the notion intersection point and critical point of $F$ interchangeably. Note also that the Reeb chords connecting $\Lambda$ to $\Lambda'$ are of two kinds: {\em short chords} which correspond to critical points of $f$ and {\em long chords} one close to each Reeb chord of $\Lambda$. As with intersection points we will sometimes identify critical points of $f$ with their corresponding short Reeb chords. Also note that to each Reeb chord of $\Lambda$ there is a unique Reeb chord of $\Lambda'$.

The following conjectural lemma is an analogue of \cite[Theorem 3.6]{EESa-duality}.
\begin{lma}[conjectural]\label{l:Morsehol}
Let $\overline{b}\,'$ denote a word of Reeb chords of $\Lambda'$ and let $\overline{b}$ denote the corresponding word of Reeb chords of $\Lambda$. Let also $\overline{e}$ denote a word of Reeb chords of $\Lambda$. If $c$ is a long Reeb chord connecting $\Lambda$ to $\Lambda'$ then let $\hat c$ denote the corresponding Reeb chord of $\Lambda$.

For sufficiently small shift $L'$ of $L$, there are the following 1-1 correspondences
\begin{itemize}
\item[$(1)$] If $a$ and $c$ are long Reeb chords then rigid disks in $\M(a;\overline{b}\,'c\,\overline{e})$ correspond to rigid disks in $\M(\hat a;\overline{b}\,\hat c\,\,\overline{e})$.
\item[$(2)$] If $a$ is a long Reeb chord and $c$ is a short Reeb chord then rigid disks in $\M(a;\overline{b}\,' c\,\overline{e})$ correspond to rigid generalized disks with disk component in $\M(\hat a;\overline{b}\,\,\overline{e})$ and with flow line asymptotic to the critical point solution of $c$ at $-\infty$ and ending at a boundary point between the last $\overline{b}$-chord and the first $\overline{e}$-chord.
\item[$(3)$] If $a$ and $c$ are short Reeb chords then rigid disks in $\M(a;c)$ correspond to rigid flow lines asymptotic to the critical point solutions of $c$ and of $a$ at $-\infty$ and $+\infty$, respectively.
\item[$(4)$] If $x$ is an intersection point and $a$ is a long Reeb chord then rigid disks in $\M(a;x)$ correspond to rigid generalized disks with disk component in $\M(\hat a)$ and with gradient line starting at $x$ and ending at the boundary of the disk.
\item[$(5)$] If $x$ is an intersection point and $c$ is a short Reeb chord then rigid disks in $\M(c;x)$ corresponds to rigid flow lines starting at $x$ and asymptotic to the critical point solution of $c$ at $+\infty$.
\end{itemize}
\end{lma}

Here $(1)-(3)$ follows from \cite[Theorem 3.6]{EESa-duality} in combination with \cite[Section 2.7]{EES-JDGgeometry} in the special case that $Y=P\times\R$ where $P$ is an exact symplectic manifold. Proofs of $(4)-(5)$ would require an analysis analogous to \cite[Section 6]{EESa-duality} carried out for a symplectization, taking into account the interpolation region used in the construction of $L'$.

\subsection{Outline of proof of Conjecture \ref{cnj:Seidel}}\label{s:pfcnj:Seidel}
Choose grading on $C(X;L,L')$ so that if $c$ is a long Reeb chord connecting $\Lambda$ to $\Lambda'$ then the degree $|c|$ satisfies $|c|=|\hat c|$ where $\hat c$ is the corresponding Reeb chord of $\Lambda$, see \eqref{e:gradingSFT}. Let $C=C_{[\alpha]}(X;L,L')$ and consider the decomposition
\[
C=C_{+}\oplus C_0,
\]
where $C_+$ is generated by long Reeb chords and where $C_0$ is generated by short Reeb chords and double points. Then $C_+$ is a sub-complex and we have the exact sequence
\[
\begin{CD}
0 @>>> C_+ @>>> C @>>> {\widehat C} @>>> 0,
\end{CD}
\]
where ${\widehat C}=C/C_+$. If $L$ and $L'$ are sufficiently close then the augmentations $\epsilon$ and $\theta$ agree and Lemma \ref{l:Morsehol} $(1)$ implies that the differential on $C_+$ is identical to that on $\V_{[\alpha]}(X,L)$. Furthermore, Lemma \ref{l:Morsehol} $(3)-(4)$ implies that the differential on ${\widehat C}$ is that of the Morse complex of $L$ and the existence of the exact sequence follows.

Consider next the isomorphism statement. Since $\C^{n}$ and $J^{1}(\R^{n-1})\times\R$ satisfy a monotonicity conditions, the homology of $C$ in a fixed degree can be computed using a fixed sufficiently large energy level. Furthermore, in $\C^{n}$ or $J^{1}(\R^{n-1})$, $L$ is displaceable, i.e. $L'$ can be moved by Hamiltonian isotopy in such a way that $L\cap L'=\emptyset$ and that there are no Reeb chords connecting $\Lambda$ and $\Lambda'$. Hence by the invariance of Lagrangian Floer cohomology proved in Theorem \ref{t:invFH}, the total complex $C$ is acyclic. The theorem follows.\qed

\subsection{Outline of proof of Corollary \ref{c:bigdiagram}}\label{s:pfc:bigdiagram}
In order to discuss Corollary \ref{c:bigdiagram} we first describe the duality exact sequence \eqref{e:seq} in more detail.

\subsubsection{Properties of the duality exact sequence}\label{ss:duality}
We recall how the exact sequence \eqref{e:seq} was constructed. Let $\Lambda''$ be a copy of $\Lambda$ shifted a large distance (compared to the length of any Reeb chord of $\Lambda$) away from $\Lambda$ in the Reeb direction and then perturbed slightly by a Morse function $f$. The part of the linearized contact homology complex of $\Lambda\cup\Lambda''$ generated by mixed Reeb chords was split as a direct sum $Q\oplus C\oplus P$. Here $Q$ is generated by the mixed Reeb chords near Reeb chords of $\Lambda$ which connects the lower sheet of $\Lambda$ to the upper sheet of $\Lambda''$, $C$ is generated by the Reeb chords which corresponds to critical points of $f$, and $P$ is generated by the Reeb chords near Reeb chords of $\Lambda$ which connects the upper sheet of $\Lambda$ to the lower sheet of $\Lambda''$. The linearized contact homology differential $\pa$ on $Q\oplus C\oplus P$ has the form
\[
\pa=
\left(\begin{matrix}
\pa_q & 0 & 0\\
\rho & \pa_c & 0\\
\eta & \sigma & \pa_p
\end{matrix}\right).
\]

Using the analogue of Lemma \ref{l:Morsehol}, the sub-complex $(P,\pa_p)$ can be shown to be isomorphic to the dual complex of the complex $(Q,\pa_q)$ using the natural pairing which pairs Reeb chords in $Q$ and $P$ which are close to the same Reeb chord of $\Lambda$. Furthermore, the complex $(Q,\pa_q)$ is canonically isomorphic to the linearized contact homology complex $(Q(\Lambda),\pa_1)$.

The next step is the observation that since $\Lambda$ is displaceable (in the sense above) the complex $C\oplus Q\oplus P$ is acyclic. Using the coarser decomposition $(Q\oplus C)\oplus P$ and writing
\[
\pa=
\left(\begin{matrix}
\pa_{qc} & 0\\
H & \pa_p
\end{matrix}\right)
\]
the chain map induced by
\[
H=\left(\begin{matrix}
\eta & \sigma
\end{matrix}\right)
\]
induces an isomorphism on homology between $Q\oplus C$ and $P$. Here the map $\sigma\colon C\to P$ counts generalized trees whereas the map from $\eta\colon Q\to P$ counts disks with boundary on $\Lambda$ and with two positive punctures disks. (To make sense of the later count and have transversely cut out moduli spaces actual disks counted have boundary on $\Lambda$ and a nearby copy $\Lambda'$.) The exact sequence \eqref{e:seq} is then constructed from the long exact sequence of the short exact sequence for the complex $C\oplus Q$.

Below we will also make use of the other splitting of the acyclic complex $Q\oplus C\oplus P$ as $Q\oplus (C\oplus P)$ with differential
\[
\pa=
\left(\begin{matrix}
\pa_{q} & 0\\
H' & \pa_{cp}
\end{matrix}\right),
\]
where
\[
H'=\left(\begin{matrix}
\rho \\ \eta
\end{matrix}\right)
\]
induces an isomorphism on homology: our proof of Corollary \ref{c:bigdiagram} relates the maps $H$ and $H'$ to isomorphisms coming from Lagrangian Floer homology.

Let $C_{\ast}$ denote the Morse complex for $f\colon \Lambda\to\R$ and let $\bar C^{\ast}$ denote the Morse complex for $-f$. Then $C_{\ast}=\bar C^{n-1-\ast}$ and we have the following diagram of chain maps
\begin{equation}\label{e:elabLCH}
\begin{CD}
@>>> \bar C^{n-k-2}@>>> P^{n-k-1} @>>> P^{n-k-1}\oplus \bar C^{n-k-1}\\
@. @| @A{H}AA @AA{H'}A\\
@>>>C_{k+1} @>>> Q_{k+1}\oplus C_{k+1} @>>> Q_{k+1}\\
@>>> \bar C^{n-k-1}\\
@. @|\,\,\cdots\\
@>>> C_k
\end{CD}
\end{equation}

\begin{lma}\label{l:commutLCH}
The diagram \eqref{e:elabLCH} commutes after passing to homology.
\end{lma}

\begin{pf}
The first and last squares commute already on the chain level by definition of the differential. Commutativity of the middle square can be seen as follows. Starting at the lower left corner with an element from $Q_{k+1}$ it is clear that the results of going up then right and right then up have common component in the $P^{n-k-1}$-summand. Using the commutativity already established we see that starting with an element in $C_{k+1}$ going up then right is the same thing as first pulling that element back to the left then going up and going two steps to the right. This vanishes in homology by exactness and hence gives the same result as going right then up. Finally, going right then up and projecting to the $\bar C^{n-k-1}$-component is the same as going right twice and then up. This vanishes in homology by exactness and gives the same as going the other way.
\end{pf}

\subsubsection{Proof of Corollary \ref{c:bigdiagram}}
We write down the diagram on the chain level. Let $(C_{\ast},\pa_{C})$ denote the Morse complex of the function $-f\colon\Lambda\to\R$ and let $(I_{\ast},\pa_{I})$ denote the Morse complex of the function $-F\colon L\to\R$ which computes the relative homology of $(\bar L,\pa\bar L)$. As above we write $\overline C^{\ast}$ and $\overline I^{\ast}$ for the corresponding co-chain complexes of $f$ and $F$, respectively.

Then the complex which computes the homology for $L$ is $(I\oplus C,\pa_{IC})$ where
\[
\pa_{IC}=
\left(
\begin{matrix}
\pa_{I} & 0\\
\mu & \pa_{C}
\end{matrix}
\right),
\]
where $\mu$ counts flow lines of $F$ starting at a critical point of $F$ and asymptotic to a critical point solution at $+\infty$. We then have the following diagram:
\begin{equation}\label{e:elabFH}
\begin{CD}
C_{k+1} @>>> I_{k+1}\oplus C_{k+1} @>>> I_{k+1} @>>> C_{k}\\
@| @| @| @|\\
\bar C^{n-k-1} @>>> \bar I^{n-k-1}\oplus \bar C^{n-k-1} @>>> \bar I _{n-k-1} @>>> \bar C_{n-k}\\
@| @V{\delta_{L,L'}}VV @VV{\delta'_{L,L'}}V @|\\
\bar C^{n-k-1}@>>> P^{n-k} @>>> P^{n-k}\oplus \bar C^{n-k} @>>> \bar C^{n-k},
\end{CD}
\end{equation}
where $\delta'_{L,L'}$ is the chain map inducing an isomorphism on homology from the splitting $C(X;L,L')=(P\oplus C)\oplus I$ instead of the splitting $P\oplus (C\oplus I)$ used in the proof sketch of Conjecture \ref{cnj:Seidel}.
Note that the bottom row of \eqref{e:elabFH} is the same as the top row of \eqref{e:elabLCH}. Joining the diagrams along this row we see that Corollary \ref{c:bigdiagram} follows once we show that \eqref{e:elabFH} commutes on the homology level. For the left and right squares this is true already on the chain level by definition of the differential. The fact that the middle square commutes for element in $\bar I^{n-k-1}$ after projection to $P^{n-k-1}$ is immediate from the definition. The same argument, using commutativity of exterior squares, as in the proof of Lemma \ref{l:commutLCH} then gives commutativity on the homology level.\qed

\end{document}